\documentclass{article}
\usepackage{tikz}
\usepackage[utf8]{inputenc}
\usepackage[T1]{fontenc}
\usepackage{mathtools}
\usepackage{hyperref}
\usepackage{stmaryrd}
\usepackage{amsmath}
\usepackage{amssymb}
\usepackage{amsfonts}
\usepackage{graphicx}
\usepackage{latexsym}
\usepackage{amsmath}
\usepackage{amssymb}
\usepackage{amsthm}
\usepackage{mathrsfs}
\usepackage{graphics}
\usepackage{xypic}
\usepackage{enumerate}
\usepackage{hyperref}
\hypersetup{
    colorlinks,
    citecolor=black,
    filecolor=black,
    linkcolor=black,
    urlcolor=black
}
\usepackage{pdfpages}
\usepackage{amssymb}
\usepackage{tikz}
\usepackage{tikz-cd}
\usepackage{appendix}
\usepackage{listings}
\usepackage{xcolor}
\usetikzlibrary{arrows.meta, positioning}
\usetikzlibrary{bending}
\usepackage{MnSymbol,bbding,pifont}
\usepackage{yfonts}
\usepackage{amsthm}
\usepackage{tikz-cd}
\usepackage{authblk}
\usepackage{verbatim}
\usepackage[a4paper, total={6in, 8in}]{geometry}
\newtheorem{theorem}{Theorem}[section]
\newtheorem{lemma}[theorem]{Lemma}
\newtheorem{corollary}[theorem]{Corollary}

\newtheorem{proposition}[theorem]{Proposition}
\newtheorem{definition}[theorem]{Definition}
\newtheorem{example}[theorem]{Example}
\newtheorem{remark}[theorem]{Remark}
\newtheorem{conjecture}[theorem]{Conjecture}

\newtheorem{problem}[theorem]{Problem}

\numberwithin{equation}{section}

\title{Diagonal orbits in the wonderful compactification}
\author{Yunsong Wei}
\date{}

\begin{document}
\setcounter{section}{-1}
\maketitle
\begin{abstract}
The various types of compactifications of symmetric spaces and locally symmetric spaces are well-studied. Among them, the De Concini-Procesi compactification, also known as the wonderful compactification, of symmetric varieties has been found to have many applications. Intuitively, this compactification provides information at infinity. The diagonal action also extends the conjugation action on semisimple groups, which has received considerable attention. In this work, we will first describe the classification of certain diagonal orbits in the wonderful compactification of a semisimple adjoint group \( G \). We will then study the compactification of the maximal torus through representations of the simply connected cover \( \tilde{G} \), which, in a sense, parameterizes these diagonal orbits. Finally, we will focus on constructing the family of closures of the Steinberg fiber. We will examine the limit of this family and show that it is a union of He-Lusztig's \( G \)-stable pieces.
\end{abstract}

\section{Introduction}People have long studied the conjugacy classes of a semisimple algebraic group \( G \). Conjugacy classes are the orbits under adjoint action, and it is natural to ask whether the adjoint orbits can be parametrized by some variety. It turns out we need to consider the adjoint invariant functions on \( G \). As a result, the algebra of invariant functions is isomorphic to the group ring of the invariant characters of a maximal torus \( T \) under the action of the Weyl group \( W \). Thus, we have the flat morphism \( G \rightarrow G/\!\!/G \cong T/W \), known as the Steinberg map. The preimage of any point in \( T/W \) is called a Steinberg fiber. The Steinberg fiber satisfies the following theorem:

\begin{theorem}[\cite{steinberg1965regular}, also summarized in \cite{humphreys1995conjugacy}]
Let \( F \) be any fiber of the Steinberg map. Then:
\begin{enumerate}
    \item \( F \) is irreducible and is a union of finitely many conjugacy classes.
    \item The semisimple elements of \( F \) form a single conjugacy class \( S \), which is closed and lies in the closure of every class in \( F \).
    \item The regular elements of \( F \) form a single conjugacy class \( R \), which is dense and open in \( F \).
    \item \( F - R \) is of codimension at least 2 in \( F \).
\end{enumerate}
\end{theorem}

Thus, the Steinberg map separates the regular orbits from the semisimple orbits, and \( T/W \) can be seen as the parameter space of adjoint orbits of \( G \). Although this is not a geometric quotient for the conjugation action, it is the best geometric model we can achieve, as Humphreys noted in \cite{humphreys1995conjugacy}.

Analogously, Kostant \cite{kostant1963lie} studied the parametrization of adjoint orbits in the semisimple Lie algebra \( \mathfrak{g} \). The algebra of \( G \)-invariant polynomial functions on \( \mathfrak{g} \) is isomorphic (under restriction) to the algebra of \( W \)-invariant polynomial functions on a Cartan subalgebra \( \mathfrak{t} \). The adjoint quotient map \( \mathfrak{g} \rightarrow \mathfrak{t} / W \) enjoys the same properties as the Steinberg map stated in the theorem above.

As Steinberg pointed out in his ICM talk \cite{steinberg1968classes}, the study of the partition of conjugacy classes in \( G \) illuminates representation theory, in addition to the algebraic and topological structure of \( G \). For example, properties invariant under conjugation, such as semisimplicity, unipotence, and nilpotence, arise from an extrinsic description through the representation of \( G \). The invariant functions of \( G \) are generated by the characters of some "basic" representation. When \( G \) is simply connected, the invariant functions are generated by the characters of its fundamental representations, which are algebraically independent, leading to \( T/W \cong \mathbb{C}^l \). These characters serve as the equations defining the Steinberg fibers.

Now consider the \( G \times G \)-equivariant compactification of \( G \) obtained by taking the closure of the image of \( G \rightarrow \mathbb{P}(\operatorname{End}(V)) \) via some representation \( V \) of \( G \). The diagonal action of \( G_{\Delta} \subset G \times G \) extends the adjoint action on \( G \) to such compactifications. One could further inquire how these orbits fit together to form the compactification.

For the wonderful compactification \( X \) of an adjoint group \( G \), the closure of a Steinberg fiber in \( X \) has been well studied (\cite{he2006unipotent}, \cite{he2006closures}, \cite{springer2006some}). As a result, every Steinberg fiber shares the same boundary in \( X \), which can be described as the disjoint union of certain \( G \)-stable pieces introduced in \cite{lusztig2004parabolic}, or, from a representation theory perspective, as elements in \( X \) that are nilpotent in any fundamental representation of \( G \). To demonstrate this, He employed a method involving the compactification of the simply connected cover \( \tilde{G} \) and determined the defining equations of the closure of the Steinberg fiber by extending the fundamental characters to the compactification \cite{he2006closures}. He then showed that the boundary corresponds exactly to the points where all the characters vanish, which is independent of the choice of the Steinberg fiber.

We aim to find similar equations that define the closure of the Steinberg fiber for adjoint type groups \( G \). A direct way of applying He's idea is described in Subsection \ref{base compactifying}, where we replace the fundamental representations with some "basic" representations of the adjoint group.

Since all the characters are pullbacks of trace functions on \( \operatorname{End}(V) \), we intend to use these trace functions to define the Steinberg fibers for the adjoint group, hoping that they will be sufficient. In invariant theory, most invariant polynomials take the form of trace functions, leading us to consider sections of the form \( \operatorname{tr}(A^i)\operatorname{tr}(Z)^j = \operatorname{tr}(A)^j\operatorname{tr}(Z^i) \), for a fixed element \( A \in \mathbb{P}(\operatorname{End}(V)) \). We expect that the common zero locus of these sections corresponds to elements \( Z \) that are 'conjugate' to \( A \) in some sense.

\noindent {\bf Lemma~\ref{conjugacy equation lemma}.} \textit{Let \( A_1, A_2 \) be two nonzero diagonal matrices in \( \operatorname{End}(\mathbb{C}^m) \) such that \( \operatorname{tr}(A_1) \neq 0 \) and \( \operatorname{tr}(A_1^e)\operatorname{tr}(A_2)^e = \operatorname{tr}(A_1)^e\operatorname{tr}(A_2^e) \) for any \( e = 1, \ldots, m \). Then there exists a constant \( t \in \mathbb{C}^* \) such that \( A_1 \) is conjugate to \( tA_2 \).}

Thus, we consider elements that are conjugate up to constants in each representation. We pose the following question:

\noindent {\bf Conjecture~\ref{Adjoint Steinberg conjecture}.} \textit{Semisimple elements of \( \tilde{G} \) are conjugate in its adjoint group \( G \) if and only if they are conjugate in \( \operatorname{GL}(V(\lambda_k)) \) up to a constant for \( k = 1, \ldots, l \).}

Here, \( V(\lambda_k) \) denotes the fundamental representations. The common zero loci of such sections consist of a single Steinberg fiber if and only if the conjecture holds. The conjecture is clearly valid for \( A_l \)-type simple groups, but requires further work to establish for other types. The conjectures for \( B_l \), \( C_l \), and \( G_2 \) cases are proven independently in Section \ref{proof of conjecture}. We heavily utilize the fact that fundamental characters separate semisimple orbits, and all the representations in these cases are self-dual.

In any case, we construct a subvariety \( Y \subset X \times \mathbb{C}^l \). General fibers of the projection \( \pi: Y \rightarrow \mathbb{C}^l \) are isomorphic to the closure of the Steinberg fiber in \( X \) when the conjecture is true. If the conjecture is false, the general fiber of \( Y \subset X \times \mathbb{C}^l \) will contain finitely many Steinberg fibers. We discuss this in Section \ref{equation}.

Since \( X \) is a compactification of \( G \), one may inquire about the limit behavior of the Steinberg fibers. A direct approach is to employ the machinery of geometric invariant theory to obtain the quotient variety \( X/\!\!/G \), as explored in \cite{he2011semistable}. It turns out that the diagonal orbits are parametrized by \( \bar{T}/W \), the orbits of the compactification of the maximal torus under the Weyl group action. The semistable points are "semisimple" points in \( X \), and the quotient map \( X \dashrightarrow \bar{T}/W \) resembles "taking the diagonal part". Another approach is to compactify the base of \( G \rightarrow T/W \) and then study the fibers, based on He's method for finding the equations of the closure of Steinberg fibers in \cite{he2011semistable}. We provide a brief discussion of these issues in Sections \ref{semistable} and \ref{base compactifying}.

We can also obtain information on the limit behavior from the projection map \( \pi: Y \rightarrow \mathbb{C}^l \). Let \( p: Y \rightarrow X \) be the projection from \( Y \subset X \times \mathbb{C}^l \) to the first factor. We show that the central fiber \( p(\pi^{-1}(0)) \) is the union of some \( G \)-stable pieces \( S_J^w \) and provide the following condition on such \( w \):

\noindent {\bf Theorem~\ref{prop9}.} \textit{The central fiber \( p(\pi^{-1}(0)) \) is the union of \( G \)-stable pieces \( S^w_J \) such that either \( \operatorname{ht}(w(\alpha_k)) > \operatorname{ht}(\alpha_k) \) or \( w(\lambda_k) \neq \lambda_k \) for each \( k \in I \backslash J \).}

We demonstrate that the central fiber cannot be pure dimensional if the rank of \( G \) exceeds 2.

\noindent {\bf Corollary~\ref{pure dimensional corollary}.} \textit{If \( G \) is simple with rank \( l \geq 3 \), then the degenerate fiber \( F_G \) is not pure dimensional.}

We compute the dimension of \( \pi^{-1}(0) \) for \( A_l \)-type groups in Section \ref{vertex}. The preimage of other subsets of \( \mathbb{C}^l \) is also discussed, for example, in Proposition \ref{preimage of subset of base}.

As \( \bar{T} \) is crucial in the study of the wonderful compactification, especially since we know that \( X/\!\!/G = \bar{T}/W \), we will study it along with the compactification of the maximal torus \( T \) through general representations of \( \tilde{G} \) in Section \ref{compactification of maximal torus}. In Section \ref{Independence of compactification}, we demonstrate that the compactification of \( T \) in \( \mathbb{P}(\operatorname{End}(V(\lambda))) \) depends only on the open face of the Weyl chamber where \( \lambda \) resides. Specifically, if the supports of \( \lambda \) and \( \lambda' \) are the same, then the compactifications of \( T \) in \( \mathbb{P}(\operatorname{End}(V(\lambda))) \) and \( \mathbb{P}(\operatorname{End}(V(\lambda'))) \) are isomorphic to each other.

\noindent {\bf Theorem~\ref{weight dependence theorem}.} \textit{The compactification of \( T \) in \( \mathbb{P}(\operatorname{End}(V(\lambda))) \) depends only on the support of \( \lambda \).}

Next, we study the normality of these toric varieties in Section \ref{normality}. In Section \ref{orbit structure of T}, we discuss the closure order of the \( T \)-orbits in its compactification within the fundamental representations of \( \tilde{G} \), as well as the irreducible components at the boundary of \( T \). As a byproduct, Theorem~\ref{boundary orbits theorem} characterizes which simple nodes correspond to the irreducible components of the boundary of \( G \) in its compactification in the fundamental representation \( \mathbb{P}(\operatorname{End}(V(\lambda_k))) \). In Section \ref{Recursive structure}, we explore the structure of \( \bar{T} \) and find that it is constructed recursively for classical groups of types \( A_l \) and \( C_l \).

In Section \ref{diagonal orbits}, we first present some general knowledge on the wonderful compactification of \( G \). We then consider the classification of diagonal orbits in \( X \backslash G \). We succeed in classifying the diagonal orbits in the minuscule \( G \times G \)-orbits, namely, those orbits corresponding to one-dimensional orbits for \( T \) acting on \( \mathbb{C}^l \).

\noindent {\bf Theorem~\ref{thm1}.} \textit{Let \( w \in W \) be the unique minimal representative for a double coset in \( W_{P_i} \backslash W / W_{P_i} \). Then the image of \( \mathfrak{g}_w \cap \mathfrak{p}_i \oplus \mathfrak{p}_i^{-} \) in \( \mathfrak{l}_i/Z(\mathfrak{l}_i) \oplus \mathfrak{l}_i/Z(\mathfrak{l}_i) = \mathfrak{sl}_2 \oplus \mathfrak{sl}_2 \) is either the diagonal or the Borel subalgebra. Moreover:
\begin{enumerate}
    \item If \( w(\alpha_i) = \alpha_i \), then the image is the diagonal.
    \item If \( w(\alpha_i) \neq \alpha_i \), then the image is the Borel subalgebra \( \mathfrak{b}^{-} \oplus \mathfrak{b} \).
\end{enumerate}}

We can use the above theorem to classify all diagonal orbits of \( X \) when \( G \) is of rank two. However, difficulties arise when considering diagonal orbits in larger \( G \times G \) orbits. In Section \ref{two simple root}, we discuss the diagonal orbits in \( G \times G \) orbits such that \( \mathfrak{l}_J \) is isomorphic to \( \mathfrak{sl}_2 \oplus \mathfrak{sl}_2 \). We see that even in such a simple case, the classification method we adopt becomes significantly more complicated.

\section*{Acknowledgement}
I would like to thank Naichung Conan Leung for sharing his ideas with me. This paper would not be finished without her insight, instructions, and encouragement. I would also like to thank Kifung Chan for helpful discussions. This work was done during the preparation of my MPhil thesis. I would like to thank the Institute of Mathematical Sciences for the marvelous research environment.
\section{Diagonal Orbits in Boundary \( G \times G \)-Orbits}\label{diagonal orbits}

\subsection{Some Recollections}
Let \( G \) be a connected semisimple group of adjoint type over the complex numbers, i.e., with a trivial center. Choose a maximal torus \( T \) of rank \( l \) and a Borel subgroup \( B \). The corresponding groups in the universal cover of \( G \) will be denoted by \( \tilde{G}, \tilde{T}, \) and \( \tilde{B} \), respectively. A Lie algebra will be denoted by lowercase letters in fraktur font corresponding to that group, giving us \( \mathfrak{t} \subset \mathfrak{b} \subset \mathfrak{g} \) respectively. Let \( X^{*}(\tilde{T}) \) be the group of characters of \( \tilde{T} \). Let \( \Phi \subset X^{*}(\tilde{T}) \) be the roots of \( \tilde{T} \) in \( \tilde{G} \), and take the positive roots \( \Phi^{+} \) to be those roots of \( \tilde{T} \) in \( \tilde{B} \) with a base of simple roots \( \Delta = \{\alpha_1, \ldots, \alpha_l\} \).

Let \( I \) be the index set \( \{1, \ldots, l\} \). The Weyl group \( W \) is generated by the simple reflections \( \{s_i\}_{i \in I} \), where \( s_i := s_{\alpha_i} \) is the reflection with respect to \( \alpha_i \). The standard parabolic subgroups correspond to the subsets of \( I \).

Let \( J \) and \( K \) be subsets of \( I \). For a subset \( J \subset I \), the Lie algebra \( \mathfrak{p}_J \) of the standard parabolic subgroup \( B \subset P_J \subset G \) has the root space decomposition
$$
\mathfrak{p}_J = \mathfrak{t} \oplus \bigoplus_{\alpha \in \Phi^+} \mathfrak{g}_{\alpha} \oplus \bigoplus_{\alpha \in \Phi^{-}(J)} \mathfrak{g}_{\beta},
$$
where \( \Phi^{-}(J) \) is the set of negative roots generated by \( \{-\alpha_i\}_{i \notin J} \).

The opposite parabolic subgroup \( P_J^{-} \) has a Lie algebra given by
$$
\mathfrak{p}_J^{-} = \mathfrak{t} \oplus \bigoplus_{\alpha \in \Phi^{+}(J)} \mathfrak{g}_{\alpha} \oplus \bigoplus_{\alpha \in \Phi^{-}} \mathfrak{g}_{\beta},
$$
where \( \Phi^{+}(J) \) is the set of positive roots generated by \( \{\alpha_i\}_{i \notin J} \). The Weyl group \( W_J \) refers to the subgroup of the Weyl group of \( \Phi \) generated by the simple reflections corresponding to \( \{\alpha_i\}_{i \notin J} \). Let \( W^J \) be the set of minimal length coset representatives of \( W / W_J \).

Let \( L_J \) be the Levi subgroup of \( P_J \) and \( G_J := L_J/Z(G_J) \). Let \( \mathfrak{u}_J \) and \( \mathfrak{u}_J^{-} \) be the nilradicals of \( \mathfrak{p}_J \) and \( \mathfrak{p}_J^{-} \), respectively, and let \( U_J, U_J^{-} \) be the corresponding subgroups of \( G \).

The support of a highest weight \( \lambda \) is defined as \( \operatorname{supp}(\lambda) = \{i \in I \mid \langle \lambda, \alpha_i \rangle \neq 0\} \). Let \( \Lambda \) be the weight lattice of \( \mathfrak{g} \), generated by the basis \( \{\lambda_1, \ldots, \lambda_l\} \) such that \( \langle \lambda_i, \alpha_j \rangle = \delta_{ij} \). The semigroup of dominant weights is denoted by \( \Lambda^{+} = \mathbb{N}\{\lambda_1, \ldots, \lambda_l\} \). If \( \lambda = \sum n_i \lambda_i \) is a dominant weight, then \( \operatorname{supp}(\lambda) = \{i \in I \mid n_i \neq 0\} \). If \( \operatorname{supp}(\lambda) = I \), then \( \lambda \) is said to be regular.

For \( w \in W \), let \( \operatorname{supp}(w) \subset I \) be the set of simple roots whose associated simple reflections occur in some (or equivalently, any) reduced decomposition of \( w \). We have the following characterization of \( \operatorname{supp}(w) \):

\begin{lemma}[\cite{he2006closures}]
Let \( w \in W \) and \( i \in I \). Then \( w \lambda_{i} \neq \lambda_{i} \) if and only if \( i \in \operatorname{supp}(w) \).
\end{lemma}

Let \( V(\lambda) \) be an irreducible \( \widetilde{G} \)-representation of dominant highest weight \( \lambda \). The representation map \( \widetilde{G} \longrightarrow \operatorname{End} V(\lambda) \backslash\{0\} \) induces a \( G \times G \)-equivariant morphism
$$
\psi_\lambda: G \longrightarrow \mathbb{P}(\operatorname{End} V(\lambda)).
$$

\begin{definition}
The wonderful compactification of \( G \) is defined as \( X = \overline{\psi_{\lambda}(G)} \subset \mathbb{P}(\operatorname{End} V(\lambda)) \) for some regular dominant highest weight \( \lambda \).
\end{definition}

For convenience, we will write \( \bar{T} \) for the closure \( \overline{\psi_{\lambda}(T)} \) when \( \lambda \) is a regular dominant highest weight in the following chapters.

\begin{theorem}[\cite{evens2008wonderful}]
The compactification is independent of the choice of the regular highest weight \( \lambda \). They are all \( G \times G \)-equivariantly isomorphic to each other.
\end{theorem}

\begin{theorem}[\cite{evens2008wonderful}]
The \( G \times G \)-orbits on the compactification correspond to subsets of \( I \). Let \( J \subset I \). Then the orbit \( S_J^\circ \), following the notation in \cite{evens2008wonderful}, is isomorphic to \( (G \times G) \times_{P_{J} \times P_{J}^{-}} G_{J} \), which \( G \times G \)-equivariantly fibers over \( G/P_J \times G/P_J^{-} \) with fiber \( G_{J} \).
\end{theorem}

Let \( z_J \) be the image of \( (1,1,1) \in (G \times G) \times_{P_{J} \times P_{J}^{-}} G_{J} \) in \( S_J^\circ \) under this isomorphism.

We will call those orbits \( S_J^\circ \) boundary orbits if \( J \neq \emptyset \) since \( S_J^\circ \subset X \backslash G \). By diagonal action, we shall always mean the action restricted to the diagonal subgroup \( G_\Delta \subset G \times G \). A number \( k \in I \) could refer to a simple reflection \( s_k \), a simple root \( \alpha_k \), the corresponding node in the Dynkin diagram of \( \mathfrak{g} \), the fundamental weight \( \lambda_k \), or the fundamental representation \( \psi_k := \psi_{\lambda_k} \), depending on the context. The Dynkin diagram of \( \mathfrak{g} \) is denoted by \( \Gamma \). All proofs for a simple algebraic group in this work can be adapted to the semisimple case.

\subsection{G-Stable Pieces}
For $J \subset I$ and $w \in W^{J}$, set $S^w_J = G_{\Delta}(Bw, 1)z_J$. Then $S_J^w$ is a locally closed subvariety of $X$ and (see \cite{he2007G})
$$
X = \coprod_{J \subset I, w \in W^{J}} S_{J}^{w}.
$$
We call $S^w_J$ a $G$-stable piece.

Let $H(\lambda)$ be the dual module for $\tilde{G}$ with lowest weight $-\lambda$. We then have a $G \times G$-equivariant morphism
$$
\rho_{\lambda}: X \rightarrow \mathbb{P}(\operatorname{End}({H}(\lambda)))
$$
which extends the morphism $G \rightarrow \mathbb{P}(\operatorname{End}({H}(\lambda)))$ defined by $g \mapsto g\left[\operatorname{Id}_{\lambda}\right]$, where $\left[\operatorname{Id}_{\lambda}\right]$ denotes the class representing the identity map on ${H}(\lambda)$ and $g$ acts by the left action.

Define ht to be the height map on the root lattice, i.e., the linear map on the root lattice which maps all the simple roots to 1.

Now fix $\lambda \in \Lambda_{+}$. Choose a basis $v_{1}, \ldots, v_{m}$ for $\mathrm{H}(\lambda)$ consisting of $T$ eigenvectors with eigenvalues $\mu_{1}, \ldots, \mu_{m}$ satisfying $\operatorname{ht}\left(\mu_{j}+\lambda\right) \geq \operatorname{ht}\left(\mu_{i}+\lambda\right)$ whenever $j \leq i$. Then $B$ is upper triangular with respect to this basis. Let $C_{J}$ be a representative of $\rho_{\lambda}\left(z_{J}\right)$ in $\operatorname{End}(\mathrm{H}(\lambda))$. Then when $\mu_{j}+\lambda$ is not a linear combination of the simple roots in $J$, we have that
$
C_{J} v_{j} \in \mathbb{C}^{\times} v_{j}.
$
If $\mu_{j}+\lambda$ is a linear combination of the simple roots and one of the simple roots is in $J$, then $C_{J} v_{j} = 0$.

\begin{proposition}
Let $\lambda = \lambda_k$ be a fundamental weight and $w \in W^{J}$. For any $b \in B$,
$
\rho_{\lambda}\left((b w, 1) z_{J}\right)
$
is represented by an upper triangular matrix.
\end{proposition}

\begin{proof}
Note that if $x$ is a nonnegative linear combination of the simple roots not in $J$, then $\operatorname{ht}(w(x)) \geqslant \operatorname{ht}(x)$. Hence,
$$
\operatorname{ht}(w(-\lambda+x)+\lambda) = \operatorname{ht}(w(x)) + \operatorname{ht}(-w \lambda+\lambda) \geqslant \operatorname{ht}(x).
$$
If the inequality is strict, then $\rho_{\lambda}\left((w, 1) z_{J}\right)$ acting on the weight vectors of weight $-\lambda+x$ is represented by an upper triangular matrix concerning the chosen basis above.

And if $\operatorname{ht}(w(-\lambda+x)+\lambda) = \operatorname{ht}(w(x)) + \operatorname{ht}(-w \lambda+\lambda) = \operatorname{ht}(x)$, then $\operatorname{ht}(w(x)) = \operatorname{ht}(x)$ and $-w \lambda + \lambda = 0$. Assume $x = a_1 \alpha_{i_1} + \cdots + a_u \alpha_{i_u}$ and $a_j \in \mathbb{N}_{>0}, i_j \in I \backslash J, j=1, \ldots, u$. Then $\operatorname{ht}(w(x)) = \operatorname{ht}(x)$ implies $w(\alpha_{i_j})$ is simple for $j=1, \ldots, u$. We want to prove $w(\alpha_{i_j}) = \alpha_{i_j}$ for $j=1, \ldots, u$. If so, let $v$ be a weight vector of weight $-\lambda+x$, then $v = x^{m_1}_{\beta_1} \cdots x^{m_s}_{\beta_s} v_m$, where $\beta_1, \ldots, \beta_s$ are positive roots which are sums of $\alpha_{i_1}, \ldots, \alpha_{i_u}$, $x_{\beta_j} \in \mathfrak{g}_{\beta_j}$ is nonzero, and $v_m$ is the lowest weight vector in $\mathrm{H}(\lambda_k)$. Then $w(x_{\beta_j}) \in \mathfrak{g}_{w(\beta_j)} = \mathfrak{g}_{\beta_j}$ implies $w(x_{\beta_j}) \in \mathbb{C} x_{\beta_j}$ as each root space is one-dimensional. And $w \lambda_k = \lambda_k$ implies $w(v_m) \in \mathbb{C} v_m$. Therefore $w(v) = w(x_{\beta_1})^{m_1} \cdots w(x_{\beta_s})^{m_s} w(v_m) \in \mathbb{C} v$.

To prove $w(\alpha_{i_j}) = \alpha_{i_j}$ for $j=1, \ldots, u$, consider the weight strings of the fundamental representation $H(\lambda_k)$; we must have $w(\alpha_k) = \alpha_k$ if $u \geq 1$. If $\alpha_{i_j}$ is connected to $\alpha_k$, then $w(\alpha_{i_j})$ is also connected to $\alpha_k$. Consider the original Dynkin diagram with $k$ removed; then the result is on a disconnected diagram. Each $s_{i}, i \neq k$ preserves the space spanned by each connected component, so $w(\alpha_{i_j})$ is a simple root in its component but connected to $\alpha_k$. Therefore, $w(\alpha_{i_j}) = \alpha_{i_j}$. The same argument applies to the simple root $\alpha_{i_{j'}}$ that is connected to $\alpha_{i_j}$. If simple root $\alpha_{i_t}$ appears in $x$, then it is connected to $\alpha_k$ and suppose the following is the unique path in the Dynkin diagram:
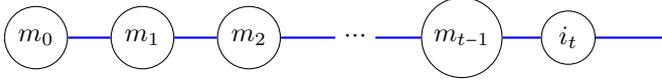
\begin{tikzpicture}[scale=.7]
  \node[circle,draw=black] (1) at (-6,0) {$m_0$};
  \node[circle,draw=black] (2) at (-4,0) {$m_1$};
  \node[circle,draw=black] (3) at (-2,0) {$m_2$};
  \node (4) at (0,0) {$\cdots$};
  \node[circle,draw=black] (5) at (2,0) {$m_{t-1}$};
  \node[circle,draw=black] (6) at (4,0) {$i_t$};
  \node (7) at (6,0) {};
  \draw [thick, blue,-] (1) -- (2) -- (3) -- (4) -- (5) -- (6) -- (7);
\end{tikzpicture}
where $k = m_0$. Then $\alpha_{m_i}$ will also appear in $x$. If $w(\alpha_{m_i}) = \alpha_{m_i}, i = 1, \ldots, m_{t-1}$, since $w(\alpha_{i_t})$ is in the component of $i_t$ in the Dynkin diagram with $k$ removed and it is the unique node connected to $m_{t-1}$ other than $m_{t-2}$, then $w(\alpha_{i_t}) = \alpha_{i_t}$. Another case is when there are three nodes $m_{t-2}, i_t$, and $i_{t+1}$ that are connected to $m_{t-1}$. This can be reduced to type $D_{t+2}(t \geq 2)$ and $k$ is the node corresponding to the standard representation of $G$. We know that there is no element in the Weyl group that exchanges the two simple roots corresponding to the two spin representations and fixes other simple roots and the highest weight of the standard representation. Finally, we get $w(\alpha_{i_j}) = \alpha_{i_j}$ for all $j = 1, \ldots, u$ by induction.

In conclusion, $\rho_{\lambda_k}\left((w, 1) z_{J}\right)(v_i)$ is either a weight vector with height larger than $v_i$ or some multiple of $v_i$, and so $\rho_{\lambda}\left((w, 1) z_{J}\right)$ is upper triangular. As a consequence, for any $b \in B$, $\rho_{\lambda_k}\left((b w, 1) z_{J}\right)$ is also represented by an upper triangular matrix.
\end{proof}

An element in $\mathbb{P}(\operatorname{End}(H(\lambda)))$ is said to be nilpotent if it may be represented by a nilpotent endomorphism of $H(\lambda)$. Define the nilpotent cone of $X$ with respect to $\lambda \in \Lambda^+$ as 
$$
\mathcal{N}(\lambda) = \left\{ z \in X \mid \rho_{\lambda}(z) \text{ is nilpotent} \right\}.
$$

We have
\begin{proposition}[\cite{he2006closures}]
$$
\mathcal{N}(\lambda) = \bigsqcup_{J \subset I} \bigsqcup_{\substack{w \in W^{J} \\ \operatorname{supp}(\lambda) \cap \operatorname{supp}(w) \neq \varnothing}} S_{J}^{w}.
$$
\end{proposition}
\begin{proof}
It is easy to see that $S_J^w \subset \mathcal{N}(\lambda)$ if and only if $w z_J \in \mathcal{N}(\lambda)$, which is equivalent to $w(\lambda) \neq \lambda$. By Lemma 10.3B of \cite{humphreys2012introduction}, $w(\lambda) = \lambda$ if and only if $w$ is generated by simple reflections that fix $\lambda$, which is equivalent to $\operatorname{supp}(w) \cap \operatorname{supp}(\lambda) = \varnothing$.
\end{proof}

Let $F \subset G$ be a Steinberg fiber, which is a fiber of the flat morphism $G \rightarrow T/W$. For example, a regular semisimple conjugacy class is a Steinberg fiber. Another example is the unipotent variety of $G$. One may refer to \cite{humphreys1995conjugacy} for general facts about the Steinberg map and Steinberg fiber.

The closure $\bar{F}$ of a Steinberg fiber in $X$ is stable under the diagonal action. We have 

\begin{theorem}[\cite{he2006closures}]
\label{thm2}
$$
\bar{F} - F = \bigcap_{k=1}^l \mathcal{N}(\lambda_k) = \bigsqcup_{J \subset I} \bigsqcup_{\substack{w \in W^{J} \\ \operatorname{supp}(w) = I}} S_{J}^{w}.
$$
\end{theorem}

We would call the elements in $\bigcap_{k=1}^l \mathcal{N}(\lambda_k)$ strongly nilpotent.

The closure relations of the $G$-stable pieces are summarized in \cite{he2007G}:
\begin{theorem}[\cite{he2007G}]
$S_J^{w_1} \subset \overline{S_K^{w_2}}$ if and only if $K \subset J$ and $w_1 \geq z^{-1} w_2 z$ for some $z \in W_{K}$.
\end{theorem}
\subsection{The Dimension Formula of Diagonal Orbits in Double Flag Varieties}

Below we give a formula for calculating the dimension of diagonal orbits of the double flag varieties $G/P_J \times G/P_J^{-}$ for later use.

Firstly, note that $G/P_J \times G/P_J^{-}$ is isomorphic to $G \times_{P_J} G/P_J^{-}$, by sending $(g_1 P_J, g_2 P_J^{-})$ to $(g_1, g_1^{-1} g_2 P_J^{-})$ and an inverse map sending $(g_3, g_4 P_J^{-})$ to $(g_3 P_J, g_3 g_4 P_J^{-})$. Thus, the orbits of the diagonal action correspond to the orbits in $G/P_J^{-}$ under the left action by $P_J$.

We have the Bruhat decomposition $B \backslash G / B^{-} = W$, and thus we have $P_J \backslash G / P_J^{-} = W_J \backslash W / W_J$. The orbit decomposition of $G/P_J \times G/P_J^{-}$ is given by
$$
G/P_J \times G/P_J^{-} = \coprod_{w \in W_J \backslash W / W_J} G \times_{P_J} P_J w P_J^{-} / P_J^{-}.
$$
Now it reduces to calculating the dimension of $P_J w P_J^{-}$ for $w \in W_J \backslash W / W_J$ as $P_J^{-}$ acts on it freely. We have 
$$
\operatorname{dim} P_J w P_J^{-} = \operatorname{dim} \mathfrak{p}_J + \operatorname{dim} w \mathfrak{p}_J^{-} w^{-1} = \operatorname{dim} \mathfrak{p}_J + \#\{\alpha \in \Phi^{-} \backslash \Phi^{-}(J) \mid w^{-1}(\alpha) \in \Phi^{-} \cup \Phi^{+}(J)\}.
$$
So we have 
\begin{proposition}
\label{prop1}
$$
\operatorname{dim} G \times_{P_J} P_J w P_J^{-} / P_J^{-} = \#\Phi^{-} \backslash \Phi^{-}(J) + \#\{\alpha \in \Phi^{-} \backslash \Phi^{-}(J) \mid w^{-1}(\alpha) \in \Phi^{-} \cup \Phi^{+}(J)\}.
$$
\end{proposition}

\subsection{Dimension of a Regular Orbit}
For an algebraic action $G$ on $X$, an orbit is said to be regular if it has maximal dimension among all orbits. An element $x \in X$ is regular if the orbit $G \cdot x$ is regular.

By the lower semi-continuity of orbit dimensions, the set of regular elements $\{x \in X \mid \operatorname{dim}(G \cdot x) \text{ is maximal}\}$ is open (see for example \cite{bate2019orbit}).

Thus, the dimension of a regular orbit in the wonderful compactification under diagonal action coincides with the dimension of a regular orbit in $G$ and equals $\operatorname{dim} G / T$.\subsection{The Classification Scheme of Esposito}
We can follow the clue of Esposito to find all diagonal orbits in the boundary divisors of the wonderful compactification of rank 2 adjoint groups or those $G \times G$-orbits $S_J^\circ$ such that $I \backslash J$ is a single element. We denote $S^\circ_{I \backslash \{i\}}$ by $S^\circ_{i}$ and call them minuscule $G \times G$-orbits.

We formulate Esposito's approach in \cite{esposito2012closures} for classifying diagonal orbits in $S_I^\circ$ through the following proposition:

\begin{proposition}[\cite{esposito2012closures}]
\label{prop2}
Consider $S_J^\circ$ as a $G \times G$-equivariant fiber bundle $\chi: (G \times G) \times_{P_{J} \times P_{J}^{-}} G_{J}$ over $G/P_J \times G/P_J^{-}$. Choose an orbit in the base space, say the orbit through $(e,w)$, then its stabilizer in $G_{\Delta}$ acts on its fiber $G_{J}$ as the image of $P_{J} \times P_{J}^{-} \cap (1, w^{-1}) G(1, w)$ in $L_J/Z(L_J) \times L_J/Z(L_J)$. All the diagonal orbits come from such a way for different $w \in W_J \backslash W / W_J$.
\end{proposition}
\subsection{Diagonal Orbits in the Minuscule $G \times G$-Orbits}

Let $P_i := P_{I \backslash \{i\}}$ be a standard minimal parabolic subgroup of the semisimple algebraic group $G$. Let $\mathfrak{p}_i$ be the minimal standard parabolic subalgebra associated with the simple roots $\alpha_i$. Let $\mathfrak{l}_i$ be the Levi subalgebra of $\mathfrak{p}_i$. Let $\mathfrak{g}^+$ denote the sum of root spaces of positive roots.

In this case, we have 
$$
\mathfrak{p}_i = \mathfrak{t} \oplus \mathfrak{g}^{+} \oplus \mathfrak{g}_{-\alpha_i}, \quad \mathfrak{l}_i = \mathfrak{t} \oplus \mathfrak{g}_{\alpha_i} \oplus \mathfrak{g}_{-\alpha_i}.
$$
Naturally, we have $\mathfrak{l}_i/Z(\mathfrak{l}_i) \oplus \mathfrak{l}_i/Z(\mathfrak{l}_i)$ is isomorphic to $\mathfrak{sl}_2 \oplus \mathfrak{sl}_2$. Denote the Lie algebra of $(1,w)G(1,w^{-1})$ by $\mathfrak{g}_w \subset \mathfrak{g} \oplus \mathfrak{g}$ for $w \in W_{P_i} \backslash W / W_{P_i}$.

\begin{lemma}
Let $w \in W$ be the unique minimal representative for a double coset in $W_{P_J} \backslash W / W_{P_J}$. Then $w(\alpha_i) \in \Phi^{+}$ and $w^{-1}(\alpha_i) \in \Phi^{+}$ for any $i \in J$.
\end{lemma}
\begin{proof}
We show this by the exchange property of the Weyl group. Say if $w(\alpha_i) = -\alpha$ or $w(\alpha) = -\alpha_i$ for some positive root $\alpha$, then $l(ws_i) < l(w)$ or $l(s_i w) < l(w)$ since $s_i$ permutes positive roots other than $\alpha_i$. By the exchange property, $w$ has a reduced expression $s_{i_1}\cdots s_{i_m}s_i$ or $s_is_{i_1}\cdots s_{i_m}$ and is not a minimal representative in $W_{P_J} \backslash W / W_{P_J}$.
\end{proof}

\begin{theorem}
\label{thm1}
Let $w \in W$ be the unique minimal representative for a double coset in $W_{P_i} \backslash W / W_{P_i}$. Then the image of $\mathfrak{g}_w \cap \mathfrak{p}_i \oplus \mathfrak{p}_i^{-}$ in $\mathfrak{l}_i/Z(\mathfrak{l}_i) \oplus \mathfrak{l}_i/Z(\mathfrak{l}_i) = \mathfrak{sl}_2 \oplus \mathfrak{sl}_2$ is either the diagonal or the Borel subalgebra. Moreover,
\begin{enumerate}
    \item If $w(\alpha_i) = \alpha_i$, then the image is the diagonal.
    \item If $w(\alpha_i) \neq \alpha_i$, then the image is the Borel subalgebra $\mathfrak{b}^{-} \oplus \mathfrak{b}$.
\end{enumerate}
\end{theorem}

\begin{proof}
For each simple root $\alpha_k$ and $X_k \in \mathfrak{g}_{\alpha_k}$, we can choose an $\mathfrak{sl}_2$-triple $X_k, Y_k, H_k$, where $Y_k \in \mathfrak{g}_{-\alpha_k}$ and $H_k = [X_k, Y_k] \in \mathfrak{t}$ such that $\alpha_k(H_k) = 2$.

We may write an element $g \in \mathfrak{p}_i$ as 
$$
g = \sum_{k=1}^{l} a_k H_k + \sum_{\alpha \in \Phi^+} u_{\alpha} X_{\alpha} + v Y_i.
$$
Here $X_{\alpha_i} = X_i$.

If $w(\alpha_i) = \alpha_i$, then $w(X_i) = X_i$, $w(Y_i) = Y_i$. And $(w(\sum_{k=1}^{l} a_k H_k) - \sum_{k=1}^{l} a_k H_k, \alpha_i) = 0$. Hence $w(\sum_{k=1}^{l} a_k H_k) - \sum_{k=1}^{l} a_k H_k \in Z(\mathfrak{l}_i)$, and the image should be diagonal.

If $w(\alpha_i) \neq \alpha_i$, then $w(\alpha_i) \in \Phi^{+} \backslash \{\alpha_i\}$. This implies $u_{\alpha_i} = 0$ since $w(X_i)$ does not lie in $\mathfrak{p}_i^{-}$.

The image of $g$ in $\mathfrak{l}_i/Z(\mathfrak{l}_i) \oplus \mathfrak{l}_i/Z(\mathfrak{l}_i)$ is thus the same as the image of $\sum_{k=1}^{l} a_k H_k + v Y_i$ as $Y_i$ will be mapped to $\mathfrak{p}_i^{-}$ by $w$ and thus survives.

The image of $w(g)$ in $\mathfrak{l}_i/Z(\mathfrak{l}_i) \oplus \mathfrak{l}_i/Z(\mathfrak{l}_i)$ is thus the same as the image of $w(\sum_{k=1}^{l} a_k H_k + u_{w^{-1}(\alpha_i)} X_{\alpha_i})$, thus the coefficient of $X_i$ survives and is independent of the coefficient $v$.

To show that $\sum_{k=1}^{l} a_k H_k$ and $w(\sum_{k=1}^{l} a_k H_k)$ have independent images in the Cartan subalgebra of $\mathfrak{sl}_2 \oplus \mathfrak{sl}_2$, we argue by contradiction. If the two are dependent, we can assume $\sum_{k=1}^{l} a_k H_k = f_i(a_1, a_2, \ldots, a_k) H_i + Z(\mathfrak{l}_i)$ and $w(\sum_{k=1}^{l} a_k H_k) = c f_i(a_1, a_2, \ldots, a_k) H_i + Z(\mathfrak{l}_i)$ for some constant $c$ and linear function $f_i$. Then $(c \sum_{k=1}^{l} a_k H_k - w(\sum_{k=1}^{l} a_k H_k), \alpha_i) = 0$, for all $a_k$, which is equivalent to $(\sum_{k=1}^{l} a_k, c \alpha_i - w^{-1}(\alpha_i)) = 0, \forall a_k, k = 1, \ldots, l$. So $w^{-1}(\alpha_i) = c \alpha_i$. But the root system is reduced, so $c = \pm 1$, a contradiction.
\end{proof}

\begin{remark}
The above discussion originated from the discussion in the paper \cite{esposito2012closures} by Esposito, where the author discusses the two kinds of images and here we show they are the only two possibilities in all ranks.
\end{remark}

Now we are ready to classify the diagonal orbits into the minuscule orbits.

Consider $S_i^\circ$ as a $G \times G$-equivariant fiber bundle $(G \times G) \times_{P_i \times P_i^{-}} G_i$ over $G/P_i \times G/P_i^{-}$. Choose an orbit in the base space, say the orbit through $(e,w)$ for $w$ being the minimal representative of $W_i \backslash W / W_i$, then the stabilizer of $(e,w) \in G/P_i \times G/P_i^{-}$ in $G_{\Delta}$ acts on its fiber $G_i$ as the image of $P_i \times P_i^{-} \cap (1, w^{-1}) G(1, w)$ in $L_i/Z(L_i) \times L_i/Z(L_i)$ by Proposition \ref{prop2}.

In this minuscule case, $G_i$ is $\mathrm{PGL_2(\mathbb{C})}$. So the above Theorem \ref{thm1} shows that if $w(\alpha_i) = \alpha_i$, then the orbits over $G_{\Delta}.(e,w)$ are divided into three classes, corresponding to the regular semisimple orbits parametrized by $\lambda \in \mathbb{C} \backslash \{0,1\} / z \sim z^{-1}$, the unipotent orbit, and the identity of dimension 2, 2, and 0 in $\mathrm{PGL}_2(\mathbb{C})$, denoted by $A_w(\lambda), U_w$, and $I_w$ respectively; if $w(\alpha_i) \neq \pm \alpha_i$, then there are exactly two orbits corresponding to the two Bruhat cells in $\mathrm{PGL}_2(\mathbb{C})$ of dimension 3 and 2, denoted by $M_w$ and $N_w$ respectively.

Therefore, $M_w$ and $N_w$ are families of Bruhat cells of $\operatorname{PGL}_2(\mathbb{C})$ over the Schubert cell $G \times_{P_i} P_i w P_i^{-} / P_i^{-}$. And $A_w(\lambda), U_w$ are families of affine quadric surfaces (i.e., $\mathbb{P}^1 \times \mathbb{P}^1 \backslash \Delta$, here $\Delta = \mathbb{P}^1$ is the diagonal of $\mathbb{P}^1 \times \mathbb{P}^1$) over the Schubert cell $G \times_{P_i} P_i w P_i^{-} / P_i^{-}$.

\begin{lemma}
The dimensions of the diagonal orbits of the minuscule double flag variety $G/P_i \times G/P_i^{-}$ are $\operatorname{dim} G \times_{P_i} P_i w P_i^{-} / P_i^{-} = \operatorname{dim} G - l - 2 - l(w)$ for $w$ being minimal representatives in $W_i \backslash W / W_i$.
\end{lemma}
\begin{proof}
Recall Proposition \ref{prop1} that 
$$
\operatorname{dim} G \times_{P_J} P_J w P_J^{-} / P_J^{-} = \#\Phi^{-} \backslash \Phi^{-}(J) + \#\{\alpha \in \Phi^{-} \backslash \Phi^{-}(J) \mid w^{-1}(\alpha) \in \Phi^{-} \cup \Phi^{+}(J)\}.
$$

If $w(\alpha_i) = \alpha_i$, then $\#\{\alpha \in \Phi^{-} \backslash \Phi^{-}(J) \mid w^{-1}(\alpha) \in \Phi^{-} \cup \Phi^{+}(J)\}$ counts those negative roots except $-\alpha_i$ which stay negative. So in total, we have $\operatorname{dim} G \times_{P_J} P_J w P_J^{-} / P_J^{-} = \operatorname{dim} G - l - l(w)$.

If $w(\alpha_i) \neq \alpha_i$ but is positive, then $\#\{\alpha \in \Phi^{-} \backslash \Phi^{-}(J) \mid w^{-1}(\alpha) \in \Phi^{-} \cup \Phi^{+}(J)\}$ still counts those negative roots except $-\alpha_i$ which stay negative since $w^{-1}(\alpha_i)$ is positive. So in total, we have $\operatorname{dim} G \times_{P_J} P_J w P_J^{-} / P_J^{-} = \operatorname{dim} G - l - l(w)$.
\end{proof}

\begin{corollary}
\label{corollary1}
The dimension formulas of the above orbits are
\begin{enumerate}
    \item 
If $w(\alpha_i) = \alpha_i$,
    \begin{description}
    \item $\operatorname{dim} A_w(\lambda) = \operatorname{dim} G - l - l(w)$,
    \item $\operatorname{dim} U_w = \operatorname{dim} G - l - l(w)$,
    \item $\operatorname{dim} I_w = \operatorname{dim} G - l - 2 - l(w)$.
    \end{description}
    \item 
If $w(\alpha_i) \neq \alpha_i$,
    \begin{description}
    \item $\operatorname{dim} M_w = \operatorname{dim} G - l + 1 - l(w)$,
    \item $\operatorname{dim} N_w = \operatorname{dim} G - l - l(w)$.
    \end{description}
\end{enumerate}
\end{corollary}
\subsection{Example: Diagonal Orbits in $G_2$}
Now we can apply the results of the previous section to classify the diagonal orbits in $G$ when $G$ is a semisimple adjoint group of rank two. Particularly, we consider $G$ of type-$G_2$.

A root system $R(G_2) \subset \mathfrak{h}^*$ associated to $G_2$ is drawn as follows:

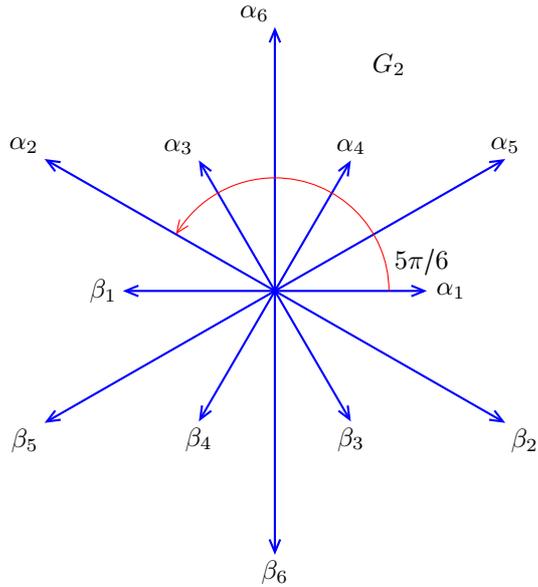
\begin{tikzpicture}[
    -{Straight Barb[bend,
       width=\the\dimexpr10\pgflinewidth\relax,
       length=\the\dimexpr12\pgflinewidth\relax]},
  ]
    \foreach \i in {0, 1, ..., 5} {
      \draw[thick, blue] (0, 0) -- (\i*60:2);
      \draw[thick, blue] (0, 0) -- (30 + \i*60:3.5);
    }
    \draw[thin, red] (1.5, 0) arc[radius=1.5, start angle=0, end angle=5*30];
    \node[right] at (2, 0) {$\alpha_1$};
    \node[above, inner sep=.2em] at (30:3.5) {$\alpha_5$};
    \node[above, inner sep=.2em] at (2*30:2) {$\alpha_4$};
    \node[above left, inner sep=.2em] at (3*30:3.5) {$\alpha_6$};
    \node[above left, inner sep=.2em] at (4*30:2) {$\alpha_3$};
    \node[above left, inner sep=.2em] at (5*30:3.5) {$\alpha_2$};
    \node[left, inner sep=.2em] at (6*30:2) {$\beta_1$};
    \node[below right, inner sep=.2em] at (11*30:3.5) {$\beta_2$};
    \node[below, inner sep=.2em] at (10*30:2) {$\beta_3$};
    \node[below, inner sep=.2em] at (9*30:3.5) {$\beta_6$};
    \node[below, inner sep=.2em] at (8*30:2) {$\beta_4$};
    \node[below left, inner sep=.2em] at (7*30:3.5) {$\beta_5$};
    \node[right] at (15:1.5) {$5\pi/6$};
    \node at (1.5, 3) {$G_2$};
\end{tikzpicture}

Where the positive roots are denoted by $\alpha_i$ and negative roots by $\beta_i = -\alpha_i$, with $\alpha_1$ and $\alpha_2$ the simple roots.

The Weyl group of $G_2$ is the dihedral group $I_2(6)$ of order 12 generated by the simple reflections $s_1$ and $s_2$ associated to $\alpha_1$ and $\alpha_2$ respectively. The Hasse digram of the Bruhat order of $I_2(6)$ is labeled as follows:

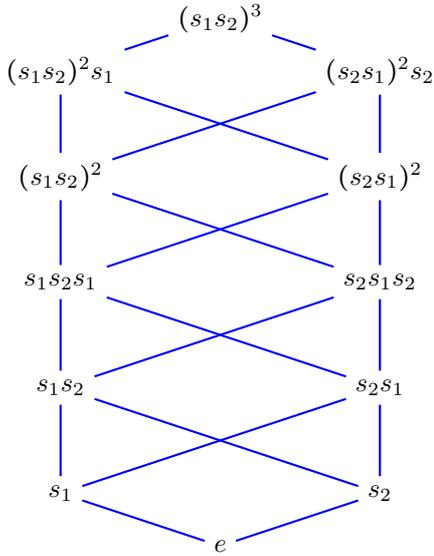
\begin{tikzpicture}[scale=.7]
  \node (top) at (0,5) {$(s_1s_2)^3$};
  \node (a) at (3,4) {$(s_2s_1)^2s_2$};
  \node (b) at (3,2) {$(s_2s_1)^2$};
  \node (c) at (3,0) {$s_2s_1s_2$};
  \node (d) at (3,-2) {$s_2s_1$};
  \node (e) at (3,-4) {$s_2$};
  \node (1) at (-3,4) {$(s_1s_2)^2s_1$};
  \node (2) at (-3,2) {$(s_1s_2)^2$};
  \node (3) at (-3,0) {$s_1s_2s_1$};
  \node (4) at (-3,-2) {$s_1s_2$};
  \node (5) at (-3,-4) {$s_1$};
  \node (bot) at (0,-5) {$e$};
  \draw [thick, blue] (bot) -- (5) -- (4) -- (3) -- (2) -- (1) -- (top) -- (a) -- (b) -- (c) -- (d) -- (e) -- (bot);
  \draw [thick, blue] (5) -- (d) -- (3) -- (b) -- (1);
  \draw [thick, blue] (e) -- (4) -- (c) -- (2) -- (a);
\end{tikzpicture}

The reflection $s_1$ sends $\alpha_1$ to $-\alpha_1$ and $\alpha_2$ to $\alpha_5 = \alpha_2 + 3\alpha_1$. The reflection $s_2$ sends $\alpha_1$ to $\alpha_3 = \alpha_1 + \alpha_2$ and $\alpha_2$ to $-\alpha_2$.

\subsubsection{Case $i=1$}
In this case, $W_{P_{1}} = \{e, s_1\}$ and the minimal representatives in $W_{P_{1}} \backslash W / W_{P_{1}}$ are $\{e, s_2, s_2s_1s_2, s_2s_1s_2s_1s_2\}$.

For $w = e$, the image of the intersection $\mathfrak{g}_w \cap \mathfrak{p}_1 \oplus \mathfrak{p}_1^{-}$ in $\mathfrak{l}_1 \oplus \mathfrak{l}_1/Z(\mathfrak{l}_1) \oplus Z(\mathfrak{l}_1)$ is the diagonal subalgebra. There is a family of regular orbits parametrized by regular adjoint orbits in $\operatorname{PSL}_2(\mathbb{C})$.

For $w = s_2$, then $w(\alpha_1) \neq \alpha_1$, so the image of the intersection $\mathfrak{g}_w \cap \mathfrak{p}_1 \oplus \mathfrak{p}_1^{-}$ in $\mathfrak{l}_1 \oplus \mathfrak{l}_1/Z(\mathfrak{l}_1) \oplus Z(\mathfrak{l}_1)$ is $\mathfrak{b}^{-} \oplus \mathfrak{b}$. And by the dimension formula in Corollary \ref{corollary1}, there is one regular orbit $G_{\Delta}.(e,s_2)z_1$.

For $w = s_2s_1s_2$, then $w(\alpha_1) \neq \alpha_1$, so the image of the intersection $\mathfrak{g}_w \cap \mathfrak{p}_2 \oplus \mathfrak{p}_2^{-}$ in $\mathfrak{l}_1 \oplus \mathfrak{l}_1/Z(\mathfrak{l}_1) \oplus Z(\mathfrak{l}_1)$ is $\mathfrak{b}^{-} \oplus \mathfrak{b}$. And by the dimension formula in Corollary \ref{corollary1}, there is no regular orbit in this case.

For $w = s_2s_1s_2s_1s_2$, then $w(\alpha_1) = \alpha_1$, so the image of the intersection $\mathfrak{g}_w \cap \mathfrak{p}_1 \oplus \mathfrak{p}_1^{-}$ in $\mathfrak{l}_1 \oplus \mathfrak{l}_1/Z(\mathfrak{l}_1) \oplus Z(\mathfrak{l}_1)$ is the diagonal subalgebra. And by the dimension formula in Corollary \ref{corollary1}, there is no regular orbit in this case.

\subsubsection{Case $i=2$}
In this case, $W_{P_{2}} = \{e, s_2\}$ and the minimal representatives in $W_{P_{2}} \backslash W / W_{P_{2}}$ are $\{e, s_1, s_1s_2s_1, s_1s_2s_1s_2s_1\}$.

For $w = e$, the image of the intersection $\mathfrak{g}_w \cap \mathfrak{p}_2 \oplus \mathfrak{p}_2^{-}$ in $\mathfrak{l}_2 \oplus \mathfrak{l}_2/Z(\mathfrak{l}_2) \oplus Z(\mathfrak{l}_2)$ is the diagonal subalgebra. There is a family of regular orbits parametrized by regular adjoint orbits in $\operatorname{PSL}_2(\mathbb{C})$.

For $w = s_1$, then $w(\alpha_2) \neq \alpha_2$, so the image of the intersection $\mathfrak{g}_w \cap \mathfrak{p}_2 \oplus \mathfrak{p}_2^{-}$ in $\mathfrak{l}_2 \oplus \mathfrak{l}_2/Z(\mathfrak{l}_2) \oplus Z(\mathfrak{l}_2)$ is $\mathfrak{b}^{-} \oplus \mathfrak{b}$. And by the dimension formula in Corollary \ref{corollary1}, there is one regular orbit $G_{\Delta}.(e,s_1)z_2$.

For $w = s_1s_2s_1$, then $w(\alpha_2) \neq \alpha_2$, so the image of the intersection $\mathfrak{g}_w \cap \mathfrak{p}_2 \oplus \mathfrak{p}_2^{-}$ in $\mathfrak{l}_2 \oplus \mathfrak{l}_2/Z(\mathfrak{l}_2) \oplus Z(\mathfrak{l}_2)$ is $\mathfrak{b}^{-} \oplus \mathfrak{b}$. And by the dimension formula in Corollary \ref{corollary1}, there is no regular orbit in this case.

For $w = s_1s_2s_1s_2s_1$, then $w(\alpha_2) = \alpha_2$, so the image of the intersection $\mathfrak{g}_w \cap \mathfrak{p}_2 \oplus \mathfrak{p}_2^{-}$ in $\mathfrak{l}_2 \oplus \mathfrak{l}_2/Z(\mathfrak{l}_2) \oplus Z(\mathfrak{l}_2)$ is the diagonal subalgebra. And by the dimension formula in Corollary \ref{corollary1}, there is no regular orbit in this case.

\subsection{Further Discussion on Diagonal Orbits in Boundary Orbits}
\subsubsection{Some General Results}
\begin{proposition}
\label{prop3}
If $w \in W_J \backslash W / W_J$ fixes all corresponding simple roots $\alpha_i, i \notin J$, then the image of $\mathfrak{p}_J \oplus \mathfrak{p}_J^{-} \cap \mathfrak{g}_w \subset \mathfrak{g} \oplus \mathfrak{g}$ in $\mathfrak{l}_{J}/Z(\mathfrak{l}_{J}) \oplus \mathfrak{l}_{J}/Z(\mathfrak{l}_{J})$ is the diagonal.
\end{proposition}

\begin{proof}
Just using the same argument in the first part of Theorem \ref{thm1}, we have $W(X_i) = X_i$ and $w(Y_i) = Y_i$ for $\alpha_i, i \notin J$, and we have $(\sum_{k=1}^{l} a_k H_k - w(\sum_{k=1}^{l} a_k H_k), \alpha_i) = 0, \forall \alpha_i, i \notin J$. Therefore, the image of $w(g)$ in $\mathfrak{l}_{J}/Z(\mathfrak{l}_{J})$ is the same as $g$.
\end{proof}

\begin{theorem}
The image of $\mathfrak{p}_J \oplus \mathfrak{p}_J^{-} \cap \mathfrak{g}_w \cap \mathfrak{t} \oplus \mathfrak{t}$ in $\mathfrak{l}_{J}/Z(\mathfrak{l}_{J}) \oplus \mathfrak{l}_{J}/Z(\mathfrak{l}_{J})$ depends only on the action of $w^{-1}$ on these simple roots $\{\alpha_i, i \notin J\}$.
\end{theorem}
\begin{proof}
Suppose $g = \sum_{k=1}^{l} a_k H_k \in \mathfrak{t}$. Assume its image in $\mathfrak{l}_{J}/Z(\mathfrak{l}_{J})$ is $\sum_{u \in I \backslash J} f_u H_u$. So $(\sum_{k=1}^{l} a_k H_k - \sum_{u \in I \backslash J} f_u H_u, \alpha_j) = 0, \forall j \in I \backslash J$. Since $\mathfrak{l}_{J}/Z(\mathfrak{l}_{J})$ is semisimple and the Cartan matrix $(\alpha_j(H_u))$ is nondegenerate, there is only one tuple of solutions independent of $w$.

Now assume the image of $w(g) = \sum_{k=1}^{l} a_k w(H_k)$ in $\mathfrak{l}_{J}/Z(\mathfrak{l}_{J})$ is $\sum_{u \in I \backslash J} g_u H_u$. Then $(\sum_{k=1}^{l} a_k w(H_k) - \sum_{u \in I \backslash J} g_u H_u, \alpha_j) = 0, \forall j \in I \backslash J$. Then $\sum_{u \in I \backslash J} g_u \alpha_j(H_u) = \sum_{k=1}^{l} a_k w^{-1}(\alpha_j)(H_k)$, and so $g_u$ depends only on the image of $w^{-1}(\alpha_j)$ for $j \in I \backslash J$.
\end{proof}

\subsubsection{Two Simple Roots Case}\label{two simple root}
We can continue to study the boundary orbits corresponding to two simple roots. It will be divided into cases based on the relation of these two simple roots so that $\mathfrak{l}_{i,j}/Z(\mathfrak{l}_{i,j})$ will be isomorphic to $\mathfrak{sl}_2 \oplus \mathfrak{sl}_2$, $\mathfrak{sl}_3$, or $\mathfrak{sp}_4$. The action of $\mathfrak{p}_J \oplus \mathfrak{p}_J^{-} \cap \mathfrak{g}_w \subset \mathfrak{g} \oplus \mathfrak{g}$ on the rank 2 semisimple Lie algebra also depends on the properties of $w$ for $w$ being minimal representatives in $W_J \backslash W / W_J$.

When $\mathfrak{l}_J \cong \mathfrak{sl}_2 \oplus \mathfrak{sl}_2$, the image of $\mathfrak{p}_J \oplus \mathfrak{p}_J^{-} \cap \mathfrak{g}_w \subset \mathfrak{g} \oplus \mathfrak{g}$ in $\mathfrak{l}_{J}/Z(\mathfrak{l}_{J}) \oplus \mathfrak{l}_{J}/Z(\mathfrak{l}_{J})$ is classified as follows:

\begin{theorem}
If $\mathfrak{l}_J \cong \mathfrak{sl}_2 \oplus \mathfrak{sl}_2$, where $J = I \backslash \{i,j\}$, and let $w \in W$ be a unique minimal representative for a double coset in $W_{P_J} \backslash W / W_{P_J}$, consider the image of $\mathfrak{p}_J \oplus \mathfrak{p}_J^{-} \cap \mathfrak{g}_w \subset \mathfrak{g} \oplus \mathfrak{g}$ in $\mathfrak{l}_{J}/Z(\mathfrak{l}_{J}) \oplus \mathfrak{l}_{J}/Z(\mathfrak{l}_{J}) = \mathfrak{sl}_2 \oplus \mathfrak{sl}_2 \oplus \mathfrak{sl}_2 \oplus \mathfrak{sl}_2$ where the 1st and 3rd entries correspond to $i$ and the 2nd and 4th entries correspond to $j$.
\begin{enumerate}
    \item If $w(\alpha_i) = \alpha_i$, $w(\alpha_j) = \alpha_j$, then the image is the sum of the diagonal in the sum of the 1st and 3rd entries and the diagonal in the sum of the 2nd and 4th entries.
    \item If $w(\alpha_i) = \alpha_i$, $w(\alpha_j) \notin \{\alpha_i, \alpha_j\}$, then the image is the sum of the diagonal in the sum of the 1st and 3rd entries and $\mathfrak{b}^- \oplus \mathfrak{b}$ in the sum of the 2nd and 4th entries.
    \item If $w(\alpha_i) = \alpha_j$, $w(\alpha_j) = \alpha_i$, then the image is the sum of the diagonal in the sum of the 1st and 4th entries and the diagonal in the 2nd and 3rd entries.
    \item If $w(\alpha_i) = \alpha_j$, and $w(\alpha_j) \notin \{\alpha_i, \alpha_j\}$, then the image is the sum of the diagonal in the sum of the 1st and 4th entries and $\mathfrak{b}^- \oplus \mathfrak{b}$ in the sum of the 2nd and 3rd entries.
    \item If $w(\alpha_i) \notin \{\alpha_i, \alpha_j\}, w(\alpha_j) \notin \{\alpha_i, \alpha_j\}$, then the image is $\mathfrak{b}^- \oplus \mathfrak{b}^- \oplus \mathfrak{b} \oplus \mathfrak{b}$ or the codimension 1 subspace of $\mathfrak{b}^- \oplus \mathfrak{b}^- \oplus \mathfrak{b} \oplus \mathfrak{b}$ containing all the root spaces and a codimension 1 Cartan subalgebra, depending on the dimension of $\langle \alpha_i, \alpha_j, w^{-1}(\alpha_i), w^{-1}(\alpha_j) \rangle$.
\end{enumerate}
\end{theorem}

\begin{proof}
The first statement follows from Proposition \ref{prop3}.

For the second, it suffices to assume $g = \sum_k a_k H_k$ since the factors in the root spaces change as we want it to be obviously. Suppose the image of $g$ is $f H_i + g H_j$ and the image of $w(g)$ is $f' H_i + g' H_j$. Then we have $(\sum_k a_k H_k - f H_i - g H_j, \alpha_i) = (\sum_k a_k H_k - f H_i - g H_j, \alpha_j) = 0 = (w(\sum_k a_k H_k) - f' H_i - g' H_j, \alpha_i) = (w(\sum_k a_k H_k) - f' H_i - g' H_j, \alpha_j)$. Since $\mathfrak{l}_J/Z(\mathfrak{l}_J) = \mathfrak{sl}_2 \oplus \mathfrak{sl}_2$, we have $\alpha_i(H_j) = \alpha_j(H_i) = 0$. And $w(\alpha_i) = \alpha_i$ implies $f' = f$. So in the 1st and 3rd entries is the diagonal. We still need to show $f, f', g, g'$ are independent as linear functions of $a_k$. So suppose $x f + y g + z g' = 0$ for any $a_k$. Then it is equivalent to that $x \alpha_i + y \alpha_j + z w^{-1}(\alpha_j) = 0$. Pairing with $\alpha_i$ we get $x = 0$ and then by the condition $w(\alpha_j) \notin \{\alpha_i, \alpha_j\}$ we know that $y = z = 0$. 

For the third, we see that the $f H_i + g H_j$ will be switched to $f H_j + g H_i$, and so do the root spaces $\mathfrak{g}_{\alpha_j}$ being switched with $\mathfrak{g}_{\alpha_i}$ and $\mathfrak{g}_{-\alpha_i}$ being switched with $\mathfrak{g}_{-\alpha_j}$.

For the fourth case, we see the root spaces go as the statement with easy checking. Assume $g = \sum_k a_k H_k$, suppose the image of $g$ is $u H_i + v H_j$ and the image of $w(g)$ is $u' H_i + v' H_j$. Then $u = \frac{1}{2} \sum_k a_k \alpha_i(H_k) = \frac{1}{2} \sum_k a_k w^{-1}(\alpha_j)(H_k) = v'$. If there exist $x, y, z$ such that $x u + y v + z u' = 0$ for any $a_k$, then $x \alpha_i + y \alpha_j + z w^{-1}(\alpha_i) = 0$. Pairing with $\alpha_i$ we get $x = 0$ and then by the condition $w^{-1}(\alpha_i) \notin \{\alpha_i, \alpha_j\}$ we know that $y = z = 0$. 

For the fifth, we still assume $g = \sum_k a_k H_k$ as root spaces move as the statement. If $\operatorname{dim} \langle \alpha_i, \alpha_j, w^{-1}(\alpha_i), w^{-1}(\alpha_j) \rangle = 4$, then $\{g(a_k), w(g)(a_k)\}$ will be the Cartan subalgebra of $\mathfrak{b}^- \oplus \mathfrak{b}^- \oplus \mathfrak{b} \oplus \mathfrak{b}$. If $\operatorname{dim} \langle \alpha_i, \alpha_j, w^{-1}(\alpha_i), w^{-1}(\alpha_j) \rangle = 3$, then $\{g(a_k), w(g)(a_k)\}$ is a codimension-one subspace. Since $\alpha_i$ and $\alpha_j$ are orthogonal, $\operatorname{dim} \langle \alpha_i, \alpha_j, w^{-1}(\alpha_i), w^{-1}(\alpha_j) \rangle \neq 2$.
\end{proof}
 \section{Compactification of Maximal Torus in Group Representations}\label{compactification of maximal torus}

\subsection{Independence of Compactification of Maximal Torus in Representations}\label{Independence of compactification}
Let $V(\lambda)$ be an irreducible representation of $\tilde{G}$ of dimension $n$ with highest weight $\lambda$ and highest weight vector $v_0$. We have the induced map:
$$\psi_{\lambda}: G \longrightarrow \mathbb{P}(\mathrm{End} V(\lambda)).$$

Choose a basis of weight vectors $\{v_i\}_{0 \leq i < n}$. Then we have a homogeneous coordinate system for $\mathbb{P}(\mathrm{End} V(\lambda))$ such that 
$$\mathbb{P}(\mathrm{End} V(\lambda)) = \{[\sum a_{ij} v_i \otimes v_j^*] \mid 0 \leq i,j < n\}.$$
We will denote the affine space of points $[\sum a_{ij} v_i \otimes v_j^*]_{0 \leq i,j < n}$ with nonzero $a_{00}$ by $\mathbb{P}_0(\mathrm{End} V(\lambda))$. Then $\psi_{\lambda}(T)$ lies in $\mathbb{P}_0(\mathrm{End} V(\lambda))$. Its closure in $\mathbb{P}_0(\mathrm{End} V(\lambda))$ is an affine toric variety, and we will denote it by $\overline{\psi_{\lambda}(T)}_0$.

Let $\mu_i$ be the weight of $v_i$ for $T$. Then $\overline{\psi_{\lambda}(T)}_0$ is isomorphic to the closure of the collection of the points 
$$(\mu_1(t)/\lambda(t), \mu_2(t)/\lambda(t), \cdots, \mu_{n-1}(t)/\lambda(t))$$ 
in $\mathbb{C}^{l-1}$. 

We also denote the set of weights of $V(\lambda)$ by $\Pi(\lambda) := \{\lambda, \mu_1, \cdots, \mu_{n-1}\}$.

We know that for $J \subset I$, the $G/P_J$ is the unique closed orbit of the action of $G$ on $\mathbb{P}V(\lambda)$ for $\lambda$ having support as $J$. Choose $p_\lambda \in \mathbb{P}V(\lambda)$ to be the point corresponding to the eigenspace with eigenvalue $\lambda$ such that $\operatorname{supp}(\lambda) = J$. Then $G.p_\lambda \cong G/P_J$. (See for example \cite{fulton2013representation}.)

Let $V(\lambda) = \oplus_{\mu \in \Pi(\lambda)} V(\lambda)_\mu$ be the weight space decomposition of $V(\lambda)$. We have the natural projections $\operatorname{pr}_{V(\lambda)_\mu}: V(\lambda) \rightarrow V(\lambda)_\mu$ to each weight space. 

\begin{lemma}
Let $J$ be a subset of $I$. Let $V(\lambda)$ be the irreducible representation with highest weight $\lambda$ with $\operatorname{supp}(\lambda) = J$. Choose $p_\lambda \in \mathbb{P}V(\lambda)$ to be the point corresponding to the eigenspace with eigenvalue $\lambda$. There is an open dense subset $O_J$ in $G/P_J$, such that for all $g \in O_{J,\lambda}$, $\operatorname{pr}_{V(\lambda)_\mu}(g.p_\lambda) \neq 0$ for all $\mu \in \Pi(\lambda)$.
\end{lemma}

\begin{proof}
For each $\mu \in \Pi(\lambda)$, there is an element $g \in G/P_J$ such that $\operatorname{pr}_{V(\lambda)_\mu}(g.p_\lambda)$ is nonzero as $V(\lambda)$ is irreducible.

Let $\{Y_\alpha\}_{\alpha \in \Phi^{-} \backslash \Phi^{-}(J)}$ be a basis of root vectors in $\mathfrak{g}/\mathfrak{p}_{J}$. Let $F: \mathbb{C}^{\dim \mathfrak{g}/\mathfrak{p}_{J}} \rightarrow G/P_J$ be the morphism sending $(t_\alpha)_{\alpha \in \Phi^{-} \backslash \Phi^{-}(J)}$ to $\operatorname{exp}(\sum_{\alpha \in \Phi^{-} \backslash \Phi^{-}(J)} t_\alpha Y_\alpha)$. Then the image of $F$ is an open dense subset of $G/P_J$. 

As $\operatorname{pr}_{V(\lambda)_\mu}(\operatorname{exp}(\sum_{\alpha \in \Phi^{-} \backslash \Phi^{-}(J)} t_\alpha Y_\alpha).p_\lambda) \neq 0$ is a polynomial open condition and for each $\mu$ there does exist $(t_\alpha)_{\alpha \in \Phi^{-} \backslash \Phi^{-}(J)}$ such that $\operatorname{pr}_{V(\lambda)_\mu}(\operatorname{exp}(\sum_{\alpha \in \Phi^{-} \backslash \Phi^{-}(J)} t_\alpha Y_\alpha).p_\lambda) \neq 0$ by choosing those close to $g \in G/P_J$ that has $\operatorname{pr}_{V(\lambda)_\mu}(g.p_\lambda) \neq 0$, we can choose 
$$O_{J,\lambda} := \{\operatorname{exp}(\sum_{\alpha \in \Phi^{-} \backslash \Phi^{-}(J)} t_\alpha Y_\alpha) \mid (t_\alpha)_{\alpha \in \Phi^{-} \backslash \Phi^{-}(J)} \in \mathbb{C}^{\dim \mathfrak{g}/\mathfrak{p}_{J}}; \operatorname{pr}_{V(\lambda)_\mu}(\operatorname{exp}(\sum_{\alpha \in \Phi^{-} \backslash \Phi^{-}(J)} t_\alpha Y_\alpha).p_\lambda) \neq 0, \forall \mu\}.$$
\end{proof}

In the complex topology, we know that the countable intersection of open dense subsets is still dense. Thus $\cap_{\lambda: \operatorname{supp}(\lambda) = J} O_{J,\lambda} \neq \emptyset$. So we have 

\begin{proposition}
Let $J$ be a subset of $I$. There is an element $g_{J}$ in $G/P_J$, such that $\operatorname{pr}_{V(\lambda)_\mu}(g.p_\lambda) \neq 0, \forall \mu \in \Pi(\lambda)$ and $\forall \lambda$ with $\operatorname{supp}(\lambda) = J$. 
\end{proposition}

Now we have $T$ acts on $G/P_J \hookrightarrow \mathbb{P}V(\lambda)$ equivariantly. Denote the closure of the orbit $T.g_{J}$ in $G/P_J$ by $T_J$. Then the closure of $T.(g_{J}.p_\lambda)$ in $\mathbb{P}V(\lambda)$ is denoted by $T_{J,\lambda}$. Then $T_J = T_{J,\lambda}$ for all $\lambda$ with $\operatorname{supp}(\lambda) = J$. Since the weights in $T_{J,\lambda}$ are the same as those in $\overline{\psi_\lambda(T)}$, we have $\overline{\psi_\lambda(T)} \cong T_{J,\lambda}$.

Thus we have proved

\begin{theorem}
\label{weight dependence theorem}
The compactification of $T$ in $\mathbb{P}(\operatorname{End}(V(\lambda)))$ depends only on the support of $\lambda$.
\end{theorem}

\begin{corollary}
If two dominant weights $\lambda$ and $\lambda'$ have the same support, then $\overline{\psi_{\lambda}(T)}_0 \cong \overline{\psi_{\lambda'}(T)}_0$.
\end{corollary}
\begin{proof}
Suppose $\lambda = \sum_{i_k \in J} n_{i_k} \lambda_{i_k}, \lambda' = \sum_{i_k \in J} m_{i_k} \lambda_{i_k}$ with $n_{i_k} > 0, m_{i_k} > 0$, and $\operatorname{supp}(\lambda) = \operatorname{supp}(\lambda') = J$. We have the $G$-equivariant embeddings
$$\prod_{i_k \in J} \mathbb{P}(V({\lambda_{i_k}})) \hookrightarrow \mathbb{P}(\otimes_{i_k \in J} V({\lambda_{i_k}})^{\otimes n_k}),$$
$$\prod_{i_k \in J} \mathbb{P}(V({\lambda_{i_k}})) \hookrightarrow \mathbb{P}(\otimes_{i_k \in J} V({\lambda_{i_k}})^{\otimes m_k}),$$
such that $(p_{\lambda_{i_k}})_{i_k \in J}$ is mapped to $\otimes_{i_k \in J} p_{\lambda_{i_k}}^{\otimes n_k}$ and $\otimes_{i_k \in J} p_{\lambda_{i_k}}^{\otimes m_k}$ respectively. The closure of the $T$-orbit of $g_J.(p_{\lambda_{i_k}})_{i_k \in J}$ will be mapped equivariantly isomorphically to $T_{J,\lambda} = \overline{\psi_{\lambda}(T)}$ and $T_{J,\lambda'} = \overline{\psi_{\lambda'}(T)}$. Thus $\overline{\psi_{\lambda}(T)}_0$ and $\overline{\psi_{\lambda'}(T)}_0$ will both correspond to those points in the closure of the $T$-orbit of $g_J.\otimes_{i_k \in J} p_{\lambda_{i_k}}^{\otimes n_k}$ such that the coordinates at $p_{\lambda_k}$ are nonvanishing for each $i_k \in J$.
\end{proof}

\begin{problem}
Prove that there is a $G \times G$-equivariant isomorphism $\overline{\psi_{\lambda}(G)} \cong \overline{\psi_{\lambda'}(G)}$ if $\operatorname{supp}(\lambda) = \operatorname{supp}(\lambda')$.
\end{problem}
\subsection{Normality of $\overline{\psi_{\lambda}(T)}_0$}
\label{normality}
Given a subset of weights $A = \{m_1, \cdots, m_s \in X^*(T)\}$, we define the map $\Phi_A: T \rightarrow \mathbb{P}^{s-1}(\mathbb{C})$ by sending $t$ to $[t^{m_1} : \cdots : t^{m_s}]$. Let $V(A) \subseteq A$ be the vertices of the convex hull of $A$ in $X^*(T) \otimes \mathbb{R}$. Then we have

\begin{proposition}
\label{prop6}
$\overline{\Phi_A(T)}$ is covered by $\overline{\Phi_A(T)} \cap \mathbb{P}^{s-1}_{k}(\mathbb{C})$, where $m_{k}$ lies in $V(A)$ and $\mathbb{P}^{s-1}_{k}(\mathbb{C}) := \{[x_1 : \cdots : x_s] \in \mathbb{P}^{s-1}(\mathbb{C}) \mid x_k \neq 0\}$.
\end{proposition}
\begin{proof}
Choose a one-parameter subgroup $n: \mathbb{C}^* \rightarrow T$, then $\Phi_A \circ n: \mathbb{C}^* \rightarrow \mathbb{P}^{s-1}(\mathbb{C})$ sends $t$ to $[t^{<m_1,n>} : \cdots : t^{<m_s,n>}].$ The limit of $t$ at infinity will be a boundary orbit $T$ at which the coordinate at $x_k$ will not vanish if $<m_k,n>$ is maximal. Since any element in the torus variety lies in the limit of some one-parameter subgroup, at some $x_k$ for $m_k \in V(A)$, the coordinate will not vanish.
\end{proof}

The proposition applies particularly to our cases of $A = \Pi(\lambda)$ for $\lambda \in \Lambda^{+}$. The vertices are $V(\Pi(\lambda))$, those highest weights. All the affine open sets corresponding to vertices of $\overline{\Phi_{\Pi(\lambda)}(T)}$ are isomorphic to $\overline{\psi_{\lambda}(T)}_0$.

So to show that $\overline{\psi_{\lambda}(T)}$ is normal, it is enough to show that $\overline{\psi_{\lambda}(T)}_0$ is normal.

\subsubsection{$A_l$-type}
In the case $G = \operatorname{PGL}_{l+1}(\mathbb{C})$, $\lambda = \lambda_k$ and $V(\lambda) = \operatorname{Sym}^r \Lambda^k \mathbb{C}^{l}$. The embedding factors through $\operatorname{Gr}(k,n) \hookrightarrow \mathbb{P}(\Lambda^k \mathbb{C}^{l}) \hookrightarrow \mathbb{P}(\operatorname{Sym}^r \Lambda^k \mathbb{C}^{l})$. For a general point $p \in \operatorname{Gr}(k,n)$ we then have $\overline{T.p}$ is isomorphic to $\overline{\psi_{\lambda_k}(T)}$ since a general point is just where those Plücker coordinates are all nonzero. 

\begin{proposition}[\cite{speyer2009matroid}]
Let $x \in G(d, n)$. Then $\overline{T x}$ is projectively normal. 
\end{proposition}

Thus, we know that $\overline{\psi_{\lambda_k}(T)} = \overline{\psi_{r\lambda_k}(T)} = \overline{T.p}$ is normal.

Another way to see this directly is as follows.

\begin{proposition}
If $\lambda_k$ is a fundamental weight of $A_l$-type, then $\overline{\psi_{\lambda_k}(T)}_0$ is a normal affine toric variety.
\end{proposition}
\begin{proof}
Consider the weight poset of $V(\lambda_k) = \Lambda^k(\mathbb{C}^{l})$, $\Pi(\lambda_k) = \{L_{i_1} + \cdots + L_{i_k} \mid 1 \leq i_1 < \cdots < i_k \leq n\}$. Then $L_{i_1} + \cdots + L_{i_k} \leq L_{j_1} + \cdots + L_{j_k}$ if and only if $i_1 \leq j_1, i_2 \leq j_2, \cdots, i_k \leq j_k$. Let 
$$A_{0} = \{L_{j_1} + \cdots + L_{j_k} - (L_{1} + \cdots + L_{k}) = (\alpha_1 + \cdots + \alpha_{j_1-1}) + (\alpha_2 + \cdots + \alpha_{j_2-1}) + \cdots + (\alpha_k + \cdots + \alpha_{j_k-1}) \mid 1 \leq j_1 \leq \cdots \leq j_k \leq n\}.$$
Then $\mathbb{N}A_0 = \{a_1\alpha_1 + \cdots + a_l\alpha_l \mid a_k \geq a_{k+1} \geq \cdots \geq a_{n} \geq 0, a_k \geq a_{k-1} \geq \cdots \geq a_{1} \geq 0\}$ is a cone in the root lattice and thus saturated.

According to Theorem 1.3.5 in \cite{cox2011toric}, we see that $\overline{\psi_{\lambda_k}(T)}_0$ is normal.
\end{proof}

\subsubsection{General Case}
\begin{conjecture}
The toric orbit closure of a general element in $G/P_{\lambda}$ is normal and can be seen from the Dynkin diagram. 

$\mathbb{N}A_0$ will still be defined by inequality and thus normal. The inequalities are like rivers centered from the supports of $\lambda$ and go along the Dynkin diagram down to the ends.
\end{conjecture}

\subsection{$T$-Orbits of $\overline{\psi_{\lambda}(T)}_0$}\label{orbit structure of T}
Let $^\lambda v_0, ^\lambda v_1, \cdots, ^\lambda v_{m_\lambda}$ be a basis of $V(\lambda)$ and the weight of $^\lambda v_{i}$ is $^\lambda \mu_i$. We could require that if $^\lambda \mu_i > ^\lambda \mu_j$, then $i < j$. Say $\lambda - ^\lambda \mu_i = \sum_{k=1}^l n_{{^\lambda \mu_i}, k} \alpha_k$ with $n_{{^\lambda \mu_i}, k} \geq 0$.

We have a natural morphism from $\mathbb{C}^l$ to $\mathbb{P}(\operatorname{End}(V(\lambda)))$ sending $(z_1, \cdots, z_l)$ to $\operatorname{diag}[1 : \prod_{k=1}^l z_k^{l_{{^\lambda \mu_1}, k}} : \cdots : \prod_{k=1}^l z_k^{l_{{^\lambda \mu_{m_\lambda}}, k}}]$. The image is exactly $\overline{\psi_{\lambda}(T)}_0$ and this morphism is $T$-equivariant. We denote this morphism by $f_\lambda: \mathbb{C}^l \rightarrow \overline{\psi_{\lambda}(T)}_0$.

The $T$-orbits of $\mathbb{C}^l$ correspond to subsets of $\{1, 2, \cdots, l\}$. The closure order of the $T$-orbits of $\mathbb{C}^l$ is isomorphic to the poset structure defined by the inclusion of subsets of $\{1, 2, \cdots, l\}$. All these data descend to the orbits of $\overline{\psi_{\lambda}(T)}_0$. Particularly, the codimension-$1$ orbits must be the image of one of the codimension-$1$ orbits in $\mathbb{C}^l$. 

Let $Z_{J}^\circ := \{(z_1, z_2, \cdots, z_l) \mid z_i \neq 0, \forall i \in I \backslash J; z_i = 0, \forall i \in J\} \subset \mathbb{C}^l$ be the orbits of $\mathbb{C}^l \backslash (\mathbb{C}^*)^l$. Let $Z_{J}$ be the closure of $Z_{J}^\circ$.

In the following, we study the number of irreducible components of $\overline{\psi_{\lambda}(T)}_0 \backslash \psi_{\lambda}(T)$ when $\lambda = \lambda_k$ is fundamental.

Firstly, number $1, 2, \cdots, l$ on the Dynkin diagrams associated with the simple roots. For simplicity, we assume the Dynkin diagram is irreducible in the following two sections.

If $j$ is an end node in the Dynkin diagram and there is a unique path connecting $k$ and $j$. Say it's of the following form:
$$
\begin{tikzpicture}[scale=.7]
  \node (1) at (-6,0) {$m_0$};
  \node (2) at (-4,0) {$m_1$};
  \node (3) at (-2,0) {$m_2$};
  \node (4) at (0,0) {$\cdots$};
  \node (5) at (2,0) {$m_{s-1}$};
  \node (6) at (4,0) {$m_{s}$};
  \draw [thick, blue] (1) -- (2) -- (3) -- (4) -- (5) -- (6);
\end{tikzpicture}$$
where $m_0 = k$, $m_s = j$.

If $\mu \in \Pi(\lambda_i)$ is a weight and $\lambda_i - \mu = \sum_{k=1}^l n_k \alpha_k$, and if the coefficients before $\alpha_{m_u}$ are not zero, then $n_{m_0}, n_{m_1}, \cdots, n_{m_{u-1}}$ are all nonzero for $u \leq s$. In fact, this statement follows from the following lemma.

\begin{lemma}[\cite{humphreys2012introduction}]
The weights in $\Pi(\lambda)$ of the form $\mu + i\alpha$ must form a connected string. If the string consists of $\mu - r\alpha, \cdots, \mu, \cdots, \mu + q\alpha$, it follows that $r - q = \langle \mu, \alpha \rangle$.
\end{lemma}

Thus, the orbits $f_{\lambda_k}(Z_{m_0}^\circ), f_{\lambda_k}(Z_{m_1}^\circ), \cdots, f_{\lambda_k}(Z_{m_{u-1}}^\circ)$ are in the boundary of $f_{\lambda_k}(Z_{u}^\circ)$. 

Let $a_\Gamma$ be the number of end nodes of the Dynkin diagram $\Gamma$ and $i_1, \cdots, i_{a_\Gamma}$ be the numbers labeled on these end nodes. From the above discussion, we see that each orbit $f_{\lambda_k}(Z_{J}^\circ)$ must lie in the boundary of some orbit $f_{\lambda_k}(Z_{\{i_s\}}^\circ)$ corresponding to the end nodes of the Dynkin diagram.

For any two orbits $f_{\lambda_k}(Z_{\{i_{s_1}\}}^\circ), f_{\lambda_k}(Z_{\{i_{s_2}\}}^\circ)$ such that $k \notin \{i_{s_1}, i_{s_2}\}$, they will not be boundary of each other since $n_{i_{s_1}}$ and $n_{i_{s_2}}$ in the expression $\lambda_{i} - \mu = \sum_{k=1}^{l} n_{k} \alpha_{k}$ are independent.

Thus we proved the following:

\begin{proposition}
\begin{enumerate}
    \item If $k = i_r$ for some $1 \leq r \leq a_\Gamma$, then the number of irreducible components in $\overline{\psi_{\lambda_k}(T)}_0 \backslash \psi_{\lambda_k}(T)$ is $a_\Gamma - 1$. The decomposition into irreducible components is 
    $$\overline{\psi_{\lambda_k}(T)}_0 \backslash \psi_{\lambda_k}(T) = f_{\lambda_k}(Z_{\{i_1\}}) \cup \cdots \cup \widehat{f_{\lambda_k}(Z_{\{i_r\}})} \cup \cdots \cup f_{\lambda_k}(Z_{\{i_{a_\Gamma}\}}).$$
    \item If $k \notin \{i_1, \cdots, i_{a_\Gamma}\}$, then the number of irreducible components in $\overline{\psi_{\lambda_k}(T)}_0 \backslash \psi_{\lambda_k}(T)$ is $a_\Gamma - 1$. The decomposition into irreducible components is 
    $$\overline{\psi_{\lambda_k}(T)}_0 \backslash \psi_{\lambda_k}(T) = f_{\lambda_k}(Z_{\{i_1\}}) \cup \cdots \cup f_{\lambda_k}(Z_{\{i_{a_\Gamma}\}}).$$
\end{enumerate}
\end{proposition}

If $k \neq i_{s}$, then it is easy to see that $f_{\lambda_k}(Z_{\{i_{s}\}}^\circ)$ does contain a $(l-1)$-dimensional torus. So $\operatorname{dim} f_{\lambda_k}(Z_{i_{s}}^\circ) = l - 1$.
\begin{example}
We number the Dynkin diagram of $A_3$ as 

\begin{tikzpicture}[scale=.7]

  \node (1) at (-6,0) {$1$};
  \node (2) at (-4,0) {$2$};
  \node (3) at (-2,0) {$3$};
 
  \draw [thick, blue] (1) -- (2) -- (3) ;
  
\end{tikzpicture}

The embedding $\psi_{\lambda_2}:\operatorname{PGL}_4(\mathbb{C})\hookrightarrow\mathbb{P}(\operatorname{End}(\bigwedge^2 \mathbb{C}^4))$ induces the map $f_{\lambda_2}:\mathbb{C}^3\rightarrow \overline{\psi_{\lambda_2}(T)}_0$ sending $(z_1,z_2,z_3)$ to $(z_2,z_2z_1,z_2z_3,z_2z_1z_3,z_2^2z_1z_3)\in\mathbb{C}^5$.
Then $f_{\lambda_2}(Z_{\{1\}})=\{(z_2,0,z_2z_3,0,0,0)|z_2,z_3\in\mathbb{C}\}$ is a $2$-dimensional irreducible component of $\overline{\psi_{\lambda_2}(T)}_0\backslash\psi_{\lambda_2}(T)$.
\end{example}

\tikzset{every picture/.style={line width=0.75pt}} 



\subsection{Irreducible Components of $\overline{\psi_{\lambda_k}(G)} \backslash \psi_{\lambda_k}(G)$}\label{boundary of G}
Now we are ready to calculate the number of irreducible components of $\overline{\psi_{\lambda_k}(G)} \backslash \psi_{\lambda_k}(G)$.

Let $X = \overline{\psi_{\lambda_0}(G)}$ for some regular dominant highest weight $\lambda_0$ be the wonderful compactification of $G$. Then we have $G \times G$-equivariant surjective morphisms $F_\lambda: X \rightarrow \overline{\psi_{\lambda}(G)}$. This morphism, when restricted to $\overline{\psi_{\lambda_0}(T)}_0 \cong \mathbb{C}^l$, is the same as $f_\lambda$ described in the last section.

According to our previous discussion on the orbits of $\overline{\psi_{\lambda_0}(T)}_0$, if $\lambda = \lambda_k$ is a fundamental weight, then the $G \times G$-orbits $F_{\lambda_k}(S_{m_{0}}^{\circ}), F_{\lambda_k}(S_{m_{1}}^{\circ}), \cdots, F_{\lambda_k}(S_{m_{u-1}}^{\circ})$ are in the boundary of $F_{\lambda_k}(S_{m_{u}}^{\circ})$. We see that each orbit $F_{\lambda_k}(S_{J}^{\circ})$ must lie in the boundary of some orbit $F_{\lambda_k}(S_{\{i_s\}}^{\circ})$ corresponding to the end nodes of the Dynkin diagram.

For a Lie algebra $\mathfrak{a}$, let $U(\mathfrak{a})$ be its enveloping algebra.

Let $i_s$ be an end node such that $k \neq i_s$. Let $V(\lambda_k)_{\{i_s\}} = U(\mathfrak{l}_{\{i_s\}}) \cdot ^{\lambda_k} v_0$. Then the following hold:

\begin{lemma}
\label{lemma1}
\begin{enumerate}
    \item $V(\lambda_k)_{\{i_s\}}$ is $\mathfrak{p}_{\{i_s\}}$-stable, and $\mathfrak{u}_{\{i_s\}}$ annihilates $V(\lambda_k)_{\{i_s\}}$.
    \item Let $Q = \{g \in G: g \cdot V(\lambda_k)_{\{i_s\}} = V(\lambda_k)_{\{i_s\}}\}$. Then $Q = P_{\{i_s\}}$.
\end{enumerate}
\end{lemma}
\begin{proof}
The proof is similar to the proof in \cite[Lemma~2.23]{evens2008wonderful}. Except a minor modification should be made on the proof of the second statement. Say there is a unique path from $k$ to $i_s$:

\begin{tikzpicture}[scale=.7]
  \node (1) at (-6,0) {$m_0$};
  \node (2) at (-4,0) {$m_1$};
  \node (3) at (-2,0) {$m_2$};
  \node (4) at (0,0) {$\cdots$};
  \node (5) at (2,0) {$m_{s-1}$};
  \node (6) at (4,0) {$i_s$};
  \draw [thick, blue] (1) -- (2) -- (3) -- (4) -- (5) -- (6);
\end{tikzpicture}
where $k = m_0$.

Then $\mathfrak{g}_{-\alpha_{i_s}}( \mathfrak{g}_{-\alpha_{m_{s-1}}} \cdots \mathfrak{g}_{-\alpha_{m_{0}}} \cdot ^{\lambda_k} v_0) \not\subset V(\lambda_k)_{\{i_s\}}$, but $\mathfrak{g}_{-\alpha_{m_{s-1}}} \cdots \mathfrak{g}_{-\alpha_{m_{0}}} \cdot ^{\lambda_k} v_0 \subset V(\lambda_k)_{\{i_s\}}$. So $Q$ is a parabolic subgroup containing $P_{\{i_s\}}$ without the root $\alpha_{i_s}$ and thus $Q = P_{\{i_s\}}$.
\end{proof}

Note that $F_{\lambda_k}(z_J) = [\operatorname{pr}_{V(\lambda_k)_{\{i_s\}}}]$. We have the following:

\begin{proposition}
The stabilizer $(G \times G)_{F_{\lambda_k}(z_{\{i_s\}})}$ of $F_{\lambda_k}(z_{\{i_s\}})$ in $G \times G$ is
$$
\left\{(xu, yv): u \in U_{\{i_s\}}, v \in U_{\{i_s\}}^{-}, x \in L_{\{i_s\}}, y \in L_{\{i_s\}}, \text{ and } xy^{-1} \in Z\left(L_{\{i_s\}}\right)\right\}.
$$
\end{proposition}
\begin{proof}
The proof is similar to the proof in \cite[Proposition~2.25]{evens2008wonderful} by using the fact proved in the above Lemma \ref{lemma1} except for the following fact that we need to verify:
$$Z\left(L_{\{i_s\}}\right) = \left\{g \in L_{\{i_s\}}: g \cdot [v] = [v], \forall [v] \in \mathbb{P}\left(V(\lambda_k)_{\{i_s\}}\right)\right\}.$$
In fact, if $g \cdot [v] = [v], \forall [v] \in \mathbb{P}\left(V(\lambda_k)_{\{i_s\}}\right)$, then $g$ is contained in $T$. For any node $u \neq i_s$ of the Dynkin diagram, there is a path starting from the node $k$ to $u$ passing through the nodes $k = n_0, n_1, \cdots, n_{s-1}$ as follows:
$$
\begin{tikzpicture}[scale=.7]
  \node (1) at (-6,0) {$n_0$};
  \node (2) at (-4,0) {$n_1$};
  \node (3) at (-2,0) {$n_2$};
  \node (4) at (0,0) {$\cdots$};
  \node (5) at (2,0) {$n_{s-1}$};
  \node (6) at (4,0) {$u$};
  \draw [thick, blue] (1) -- (2) -- (3) -- (4) -- (5) -- (6);
\end{tikzpicture}
$$
None of the $n_0, \cdots, n_{s-1}$ equals $i_s$. Then we have $\alpha_{n_0}(g) = \alpha_{n_0}(g)\alpha_{n_1}(g) = \cdots = \alpha_{n_0}(g)\alpha_{n_1}(g) \cdots \alpha_{n_{s-1}}(g)\alpha_{u}(g) = 1$. This is equivalent to $\alpha_{u}(g) = 1$ for any $u \neq i_s$, which means exactly that $g \in Z\left(L_{\{i_s\}}\right)$.
\end{proof}

\begin{corollary}
If $k \neq i_s$, then $F_{\lambda_k}$ induces a $G \times G$-equivariant isomorphism between $F_{\lambda_k}(S_{\{i_{s}\}}^\circ)$ and $S_{\{i_{s}\}}^\circ$.
\end{corollary}

Particularly, we have $\operatorname{dim} F_{\lambda_k}(S_{\{i_{s}\}}^\circ) = \operatorname{dim} G - 1$ and thus we obtained the following theorem:

\begin{theorem}\label{boundary orbits theorem}
\begin{enumerate}
    \item If $k = i_r$ for some $1 \leq r \leq a_\Gamma$, then the number of irreducible components in $\overline{\psi_{\lambda_k}(G)} \backslash \psi_{\lambda_k}(G)$ is $a_\Gamma - 1$. The decomposition into irreducible components is 
    $$\overline{\psi_{\lambda_k}(G)} \backslash \psi_{\lambda_k}(G) = F_{\lambda_k}(S_{\{i_1\}}) \cup \cdots \cup \widehat{F_{\lambda_k}(S_{\{i_r\}})} \cup \cdots \cup F_{\lambda_k}(S_{\{i_{a_\Gamma}\}}).$$
    \item If $k \notin \{i_1, \cdots, i_{a_\Gamma}\}$, then the number of irreducible components in $\overline{\psi_{\lambda_k}(G)} \backslash \psi_{\lambda_k}(G)$ is $a_\Gamma - 1$. The decomposition into irreducible components is 
    $$\overline{\psi_{\lambda_k}(G)} \backslash \psi_{\lambda_k}(G) = F_{\lambda_k}(S_{\{i_1\}}) \cup \cdots \cup F_{\lambda_k}(S_{\{i_{a_\Gamma}\}}).$$
\end{enumerate}
\end{theorem}
\begin{proof}
We prove the first case. Due to the previous discussion, we just need to show $F_{\lambda_k}(S_{\{i_{m_1}\}}^\circ) \cap F_{\lambda_k}(S_{\{i_{m_2}\}}^\circ) = \varnothing$.

If $1 \leq m_1 \leq m_2 \leq i_{a_\Gamma}$ and $i_r \notin \{m_1, m_2\}$ if $k = i_r$ for some $1 \leq r \leq a_\Gamma$. Since $S_{\{i_{m_1}\}}^\circ = \coprod_{w \in W^{\{i_{m_1}\}}} G.(Bw, 1) z_{\{i_{m_1}\}}$, we have 
$$F_{\lambda_k}(S_{\{i_{m_1}\}}^\circ) \cap \overline{\psi_{\lambda_k}(T)_{0}} = F_{\lambda_k}(\coprod_{w \in W^{\{i_{m_1}\}}, w(\alpha_i) = \alpha_i, \forall i \neq i_{m_1}} (Tw, 1) z_{\{i_{m_1}\}})$$
And
$$F_{\lambda_k}(\coprod_{w \in W^{\{i_{m_1}\}}, w(\alpha_i) = \alpha_i, \forall i \neq i_{m_1}} (Tw, 1) z_{\{i_{m_1}\}}) = F_{\lambda_k}((T, 1) z_{\{i_{m_1}\}}).$$
As $F_{\lambda_k}((T, 1) z_{\{i_{m_1}\}}) \neq F_{\lambda_k}((T, 1) z_{\{i_{m_2}\}})$ in the toric variety $\overline{\psi_{\lambda_k}(T)_{0}}$, we have $F_{\lambda_k}(S_{\{i_{m_1}\}}^\circ) \cap F_{\lambda_k}(S_{\{i_{m_2}\}}^\circ) = \varnothing$.
\end{proof}

\begin{example}
Let $G = \operatorname{PGL}_{l+1}(\mathbb{C})$. The fundamental representations are $V(\lambda_k) = \bigwedge^k \mathbb{C}^{l+1}, 1 \leq k \leq l$. Then there are two irreducible components in $\overline{\psi_{\lambda_k}(\operatorname{PGL}_{l+1}(\mathbb{C}))} \backslash \psi_{\lambda_k}(\operatorname{PGL}_{l+1}(\mathbb{C}))$ if $1 < k < l$. The components also come from the following two natural rational maps:
$$\mathbb{P}(\operatorname{End}(\mathbb{C}^{l+1})) \dashrightarrow \mathbb{P}(\operatorname{End}(\bigwedge^k \mathbb{C}^{l+1}))$$
$$\mathbb{P}(\operatorname{End}(\bigwedge^{l} \mathbb{C}^{l+1})) \dashrightarrow \mathbb{P}(\operatorname{End}(\bigwedge^{l+1-k} \bigwedge^{l} \mathbb{C}^{l+1})).$$

Let $V_1 = \langle e_1 \wedge e_{i_1} \wedge \cdots \wedge e_{i_{k-1}} \rangle_{2 \leq i_1 < i_2 < \cdots < i_{k-1} \leq l+1}$, $V_2 = \langle e_{j_1} \wedge \cdots \wedge e_{j_{k}} \rangle_{2 \leq j_1 < j_2 < \cdots < j_{k} \leq l+1}$. Then the two components correspond to $G \times G.[\operatorname{pr}_{V_1}]$ and $G \times G.[\operatorname{pr}_{V_2}]$.

A way to see that $V_1$ and $V_2$ are not in one $G \times G$-orbit is the following: Suppose $r[\operatorname{pr}_{V_1}] s = [\operatorname{pr}_{V_2}]$ for some $r, s \in G$. Then $s(V_2) \subset V_1$ and $s|_{V_2}: V_2 \rightarrow V_1$ is an isomorphism. Under the basis $\{e_1, \cdots, e_{l+1}\}$, we have minors $\operatorname{det}_{h_1, \cdots, h_k}^{j_1, \cdots, j_k}(g) = 0$ for all $2 \leq h_1 < \cdots < h_k \leq l+1, 2 \leq j_1 < \cdots < j_k \leq l+1$. But this implies $\operatorname{rank}(g) \leq k + 1 < l + 1$.
\end{example}

\begin{problem}
Find the pattern for the number of irreducible components in $\overline{\psi_{\lambda}(G)} \backslash \psi_{\lambda}(G)$ for general highest dominant weight $\lambda$.
\end{problem}
\subsection{Recursive Structure of the Toric Variety}\label{Recursive structure}
Given a root system, the fan generated by the Weyl chambers is associated with a smooth complete toric variety. Thus, we have a functor $M$ from the collection of root systems to the collection of toric varieties \cite{batyrev2011functor}. For example, when the root system is $A_1$, the toric variety $M(A_1)$ is $\mathbb{P}^1$. Recall that in the wonderful compactification, the closure of $\psi_{\lambda}(\mathrm{T})$ for $\lambda$ regular is the toric variety associated with the root system of $\mathfrak{g}$. We will see that they are successive blow-ups from projective space in type $A_l$, and they are constructed in a recursive way for the two families of adjoint-type classic groups, $A_l$, and $C_l$.

\subsubsection{$A_l$}
Adding a factor of projective space of endomorphisms of fundamental representations will yield more blow-ups. To see this, let $\Psi_i$ be the embedding $\mathrm{PGL}_{l+1}(\mathbb{C}) \hookrightarrow \mathbb{P}(\mathrm{End}(\Lambda^1(\mathbb{C}^{l+1}))) \times \cdots \times \mathbb{P}(\mathrm{End}(\Lambda^i(\mathbb{C}^{l+1})))$ and $\pi_i$ be the canonical projection 
$$ \mathbb{P}(\mathrm{End}(\Lambda^1(\mathbb{C}^{l+1}))) \times \cdots \times \mathbb{P}(\mathrm{End}(\Lambda^i(\mathbb{C}^{l+1}))) \rightarrow \mathbb{P}(\mathrm{End}(\Lambda^1(\mathbb{C}^{l+1}))) \times \cdots \times \mathbb{P}(\mathrm{End}(\Lambda^{i-1}(\mathbb{C}^{l+1}))) $$ 
for $i \leq n$. Then we have 

\begin{proposition}
$\overline{\Psi_2(\mathrm{T})}$ is the blow-up of $\overline{\Psi_1(\mathrm{T})} = \mathbb{P}^{l}$ along the $l+1$ points fixed by $\mathrm{T}$ in $\mathbb{P}^{l}$. And inductively, $\overline{\Psi_{i+1}(\mathrm{T})}$ is the blow-up of $\overline{\Psi_{i}(\mathrm{T})}$ along the proper transform of those $(i-1)$-dimensional hyperplanes in $\mathbb{P}^{l}$ fixed by some codimension-$i-1$ subtori of $\mathrm{T}$. 
\end{proposition}

Particularly, $\overline{\Psi_{l}(\mathrm{T})} = X(A_l)$ is obtained by first blowing up all the fixed points, then the 1-dimensional lines connecting two of the fixed points, and then 2-dimensional planes containing any three of the fixed points.

Besides, we have an inductive structure of $X(A_*)$.

\begin{proposition}
Consider the projection $\pi_1 \circ \cdots \circ \pi_{l}: X(A_l) \rightarrow \mathbb{P}^{l}$. The preimage of a point in the $(i-1)$-dimensional hyperplane fixed by some codimension-$i-1$ subtori but not fixed by any codimension-$i$ subtori in $\mathbb{P}^{l}$ is isomorphic to $X(A_{l-i})$.
\end{proposition}

The proofs are by induction, so instead, we would like to demonstrate an example here.

\begin{example}
Consider the case $n=3$. $(t_1, t_2, t_3, t_4) \in \mathrm{T}$ is mapped to 
$$([t_1, t_2, t_3, t_4], [t_1 t_2, t_1 t_3, t_1 t_4, t_2 t_3, t_2 t_4, t_3 t_4], [t_1 t_2 t_3, t_1 t_2 t_4, t_1 t_3 t_4, t_2 t_3 t_4])$$ 
by $\Psi_3$, so the preimage of $(1, t_2, t_3, t_4)$ contains 
$$([1, t_2, t_3, t_4], [t_2, t_3, t_4, t_2 t_3, t_2 t_4, t_3 t_4], [t_2 t_3, t_2 t_4, t_3 t_4, t_2 t_3 t_4])$$ 
which has the same vanishing order as $([t_2, t_3, t_4], [t_2 t_3, t_2 t_4, t_3 t_4])$ and whose closure is exactly $X(A_2)$. 

And the preimage of $([1, 2, t_3, t_4])$ contains 
$$([1, 2, t_3, t_4], [2, t_3, t_4, 2 t_3, 2 t_4, t_3 t_4], [2 t_3, 2 t_4, t_3 t_4, 2 t_3 t_4])$$ 
which has the same vanishing order as $([t_3, t_4])$ and is exactly $X(A_1)$.
\end{example}

\subsubsection{Classic Groups of Types $A_l$ and $C_l$}
Let $\Psi_i^{l}$ be the embedding for one of the following cases:

\begin{enumerate}
    \item $\mathrm{PSL}_{l+1}(\mathbb{C}) \hookrightarrow \mathbb{P}(\mathrm{End}(\Lambda^1(\mathbb{C}^{l+1}))) \times \cdots \times \mathbb{P}(\mathrm{End}(\Lambda^i(\mathbb{C}^{l+1})))$ for $A_{l}$-type
    \item $\mathrm{PSp}_{2l}(\mathbb{C}) \hookrightarrow \mathbb{P}(\mathrm{End}(\Lambda^1(\mathbb{C}^{2l}))) \times \cdots \times \mathbb{P}(\mathrm{End}(\Lambda^i(\mathbb{C}^{2l})))$ for $C_l$-type
\end{enumerate}

Let $\pi_i^{l}$ be the canonical projection for one of the following cases respectively:

\begin{enumerate}
    \item $\mathbb{P}(\mathrm{End}(\Lambda^1(\mathbb{C}^{l+1}))) \times \cdots \times \mathbb{P}(\mathrm{End}(\Lambda^i(\mathbb{C}^{l+1}))) \rightarrow \mathbb{P}(\mathrm{End}(\Lambda^1(\mathbb{C}^{l+1}))) \times \cdots \times \mathbb{P}(\mathrm{End}(\Lambda^{i-1}(\mathbb{C}^{l+1})))$ for $A_l$-type
    \item $\mathbb{P}(\mathrm{End}(\Lambda^1(\mathbb{C}^{2l}))) \times \cdots \times \mathbb{P}(\mathrm{End}(\Lambda^i(\mathbb{C}^{2l}))) \rightarrow \mathbb{P}(\mathrm{End}(\Lambda^1(\mathbb{C}^{2l}))) \times \cdots \times \mathbb{P}(\mathrm{End}(\Lambda^{i-1}(\mathbb{C}^{2l})))$ for $C_l$-type
\end{enumerate}

Let $\mathrm{T}$ be the subgroup of diagonal matrices of $\mathrm{PSL}_{l+1}(\mathbb{C})$ or $\mathrm{PSp}_{2l}(\mathbb{C})$ respectively. Then we have 

\begin{proposition}
The preimage of a fixed point of $\overline{\Psi_1^{l}(\mathrm{T})}$ under $\pi_2^{l}$ is isomorphic to $\overline{\Psi_1^{l-1}(\mathrm{T})}$.
\end{proposition}

This is a particular case of the following statement:

\begin{proposition}
The preimage of a fixed point of $\overline{\Psi_1^{l}(\mathrm{T})}$ under $\pi_2^{l} \circ \pi_3^{l} \circ \cdots \circ \pi_i^{l}$ is isomorphic to $\overline{\Psi_{i-1}^{l-1}(\mathrm{T})}$.
\end{proposition}

More generally:

\begin{proposition}
The preimage of a point of $\overline{\Psi_1^{l}(\mathrm{T})}$ in the locus fixed by some codimension-$j$ subtori but not fixed by any codimension-{(j-1)} subtori under $\pi_2^{l} \circ \pi_3^{l} \circ \cdots \pi_i^{l}$ is isomorphic to $\overline{\Psi_{i-j-1}^{l-j-1}(\mathrm{T})}$.
\end{proposition}

Another generalization:

\begin{proposition}
The preimage of a fixed point of $\overline{\Psi_l^{l}(\mathrm{T})}$ under $\pi_{l+1}^{l} \circ \pi_{l+2}^{l} \circ \cdots \pi_i^{l}$ is isomorphic to $\overline{\Psi_{i-l}^{l-l}(\mathrm{T})}$.
\end{proposition}

In particular:

\begin{corollary}
The preimage of a point of $\overline{\Psi_1^{l}(\mathrm{T})}$ in the locus fixed by some codimension-$j$ subtori but not fixed by any codimension-{(j-1)} subtori under $\pi_2^{l} \circ \pi_3^{l} \circ \cdots \pi_{l}^{l}$ is isomorphic to $X(A_{l-j-1})$ or $X(C_{l-j-1})$ respectively.
\end{corollary}

Thus, the toric varieties $X(A_{l})$ and $X(C_{l})$ are determined by their different "initial conditions" as $\overline{\Psi_1^{l}(\mathrm{T})}$ with the same recursive construction.

\section{Family of Diagonal Orbits in Wonderful Compactifications}\label{family}

\subsection{Equations Defining Steinberg Fibers}\label{equation}
The wonderful compactification can also be defined as the closure of the image of 
$$\rho: G \rightarrow \mathbb{P}(\operatorname{End}(H(\lambda_1))) \times \cdots \times \mathbb{P}(\operatorname{End}(H(\lambda_l))).$$
Here, $\rho$ can be viewed as $(\rho_{\lambda_1}, \cdots, \rho_{\lambda_l})$.

Choose a basis $^k v_{1}, \ldots, ^k v_{m_k}$ for $\mathrm{H}(\lambda_k)$ consisting of $T$-eigenvectors with eigenvalues $^k \mu_{1}, \ldots, ^k \mu_{m_k} = -\lambda_k$, satisfying $\mathrm{ht}(^k \mu_{j} + \lambda_k) \geq \operatorname{ht}(^k \mu_{i} + \lambda_k)$ whenever $j \leq i$. Then $^k \mu_j + \lambda_k$ is the nonnegative sum of simple roots. 

We have an open embedding $a: \mathbb{C}^l \hookrightarrow \mathbb{P}(\operatorname{End}(H(\lambda_1))) \times \cdots \times \mathbb{P}(\operatorname{End}(H(\lambda_l)))$, which factors through $X$, sending $(z_1, \cdots, z_l) \in \mathbb{C}^l$ to 
$$(\operatorname{diag}[\prod_{u=1}^{l} z_u^{\langle \lambda_u, ^1 \mu_1 + \lambda_1 \rangle} : \cdots : \prod_{u=1}^{l} z_u^{\langle \lambda_u, ^1 \mu_{m_1} + \lambda_1 \rangle}], \cdots, \operatorname{diag}[\prod_{u=1}^{l} z_u^{\langle \lambda_u, ^l \mu_1 + \lambda_l \rangle} : \cdots : \prod_{u=1}^{l} z_u^{\langle \lambda_u, ^l \mu_{m_l} + \lambda_l \rangle}])$$

We construct the following family: 
$$X_z := \{g \in X \mid f_{z, e_1, e_2, k}(g) = 0,\, 1 \leq e_1 < e_2 \leq \operatorname{dim}(H(\lambda_k)),\, 1 \leq k \leq l\}$$
where 
$$f_{z, e_1, e_2, k}(g) := \left(\sum_{1 \leq j \leq m_k} \left(\prod_{u=1}^{l} z_u^{\langle \lambda_u, ^k \mu_j + \lambda_k \rangle}\right)^{e_1}\right)^{\operatorname{l.c.m}(e_1, e_2)/e_1} \operatorname{tr}_{H(\lambda_k)}(\rho_{\lambda_k}(g)^{e_2})^{\operatorname{l.c.m}(e_1, e_2)/e_2}$$

$$ - \left(\sum_{1 \leq j \leq m_k} \left(\prod_{u=1}^{l} z_u^{\langle \lambda_u, ^k \mu_j + \lambda_k \rangle}\right)^{e_2}\right)^{\operatorname{l.c.m}(e_1, e_2)/e_2} \operatorname{tr}_{H(\lambda_k)}(\rho_{\lambda_k}(g)^{e_1})^{\operatorname{l.c.m}(e_1, e_2)/e_1}$$
and $z = (z_1, \cdots, z_l) \in \mathbb{C}^l$.

Then $f_{z, e_1, e_2, k}(g) = 0,\, 1 \leq e_1 < e_2 \leq \operatorname{dim}(H(\lambda_k))$ if and only if $\rho_{\lambda_k}(g)$ is nilpotent in $\mathbb{P}(\operatorname{End}(H(\lambda_k)))$ or the semisimple part of $\rho_{\lambda_k}(g)$ is conjugate to $\operatorname{diag}[\prod_{u=1}^{l} z_u^{\langle \lambda_u, ^k \mu_1 + \lambda_k \rangle} : \cdots : \prod_{u=1}^{l} z_u^{\langle \lambda_u, ^k \mu_{m_k} + \lambda_k \rangle}]$ in $\mathbb{P}(\operatorname{End}(H(\lambda_k)))$.

\begin{proposition}
\label{prop5}
If $G = \operatorname{PGL}_{l+1}(\mathbb{C})$, then the closure of the diagonal orbit containing $a(z)$ for general $z \in \mathbb{C}^{l}$ is $X_z$.
\end{proposition}
\begin{proof}
It is easy to see that the diagonal orbits in $\mathbb{P}(\operatorname{End}(H(\lambda_l))) = \mathbb{P}(\operatorname{End}(\mathbb{C}^{l+1}))$ satisfying the equations $f_{z, e_1, e_2, l} = 0$ for all $ 1 \leq e_1 < e_2 \leq \operatorname{dim}(H(\lambda_k))$ contain exactly one regular semisimple orbit, which is the orbit through $a(z)$. Thus, if $X_z$ contains other regular orbits of $X$, then the regular orbits must be in the boundary orbits $X \backslash G$ and cannot be of full rank as projective transformations of each representation. Therefore, they must be nilpotent on each fundamental representation. However, we know such orbits are boundaries of any regular semisimple orbits and therefore cannot be regular (Theorem \ref{thm2}).
\end{proof}

\begin{proposition}
If $G = \operatorname{PSp}_{2l}(\mathbb{C})$ or $\operatorname{SO}_{2l+1}(\mathbb{C})$, then the closure of the diagonal orbit containing $a(z)$ for general $z \in \mathbb{C}^l$ is $X_z$.
\end{proposition}
\begin{proof}
According to the proof of Proposition \ref{prop5}, it suffices to find a $k$ and $\mathbb{P}(\operatorname{End}(H(\lambda_k)))$ such that the variety defined by the equations $f_{z, e_1, e_2, l} = 0$ for all $ 1 \leq e_1 < e_2 \leq \operatorname{dim}(H(\lambda_k))$ contains only one regular semisimple orbit. 

Let $u = 2l$ if $G = \operatorname{PSp}_{2l}(\mathbb{C})$, and $u = 2l + 1$ if $G = \operatorname{SO}_{2l+1}(\mathbb{C})$.

We may choose $k$ such that $H(\lambda_k) = \mathbb{C}^{u}$ is the standard representation of the classical group $G$. Let $\rho_k: G \rightarrow \mathbb{P}(\operatorname{End}(\mathbb{C}^{u}))$ be the embedding. 

For general $t, s \in T$, if $\rho_k(t)$ and $\rho_k(s)$ are in the same orbit under the natural action of $S_{u}$ in $\mathbb{P}(\operatorname{End}(\mathbb{C}^{u}))$, then there is a group element $w \in W$ such that $w(t) = s$. Hence, they are in the same regular semisimple orbit. For example, it is evident that if $G = \operatorname{PSp}_{2l}(\mathbb{C})$, then $t = [t_1 : t_2 : t_3 : t_4]$ is mapped to $[t_1, t_2, t_3, t_4]$ with the same set of eigenvalues (up to scalar) if and only if they are conjugate under the Weyl group of $\operatorname{PSp}_{2l}(\mathbb{C})$. 
\end{proof}

The above method does not work for $\operatorname{PSO}_{2l}(\mathbb{C})$ since its Weyl group will contain elements acting on each coordinate independently. For example, $[2:3:4:1/2:1/3:1/4]$ is not conjugate to $[1/2:3:4:2:1/3:1/4]$ under the conjugate action of $\operatorname{PSO}_{4}(\mathbb{C})$ in $\mathbb{P}(\operatorname{End}(\mathbb{C}^4))$.

Recall the Steinberg map $\operatorname{St}: G \rightarrow T/W$ defined by $x \mapsto [x_s]$, which sends a group element to the conjugacy classes of its semisimple component. This map is induced from the isomorphism $\mathbb{C}[G]^G \cong \mathbb{C}[T]^W$. A Steinberg fiber is any fiber of $\operatorname{St}$. For the simply connected case, $\mathbb{C}[\tilde{G}]^{\tilde{G}}$ is freely generated by the characters of $\tilde{G}$'s fundamental representations, so that any two group elements lie in the same Steinberg fiber if and only if they have the same trace in each fundamental representation. As the equations defining $X_z$ are related to these characters, we hope that $X_z$ is the closure of some Steinberg fiber of $G$ in $X$ for $z \in (\mathbb{C}^*)^l$.

Fix $z \in (\mathbb{C}^*)^l$. Let $g \in G$ and $g = g_s g_u$ be its Jordan decomposition. Then $g \in X_z$ if and only if $g_s \in X_z$. So we have that $X_z$ is the closure of a Steinberg fiber if and only if for any two $t_1, t_2 \in T$, $t_1, t_2 \in X_z$ implies $t_1 \in W.t_2$.

However, the equations defining $X_z$ only control the set of eigenvalues up to a constant in each representation if $z$ is general enough. 

For example, if $z$ satisfies $\sum_{1 \leq j \leq m_k} \left( \prod_{u=1}^{l} z_u^{\langle \lambda_u, ^k \mu_j + \lambda_k \rangle} \right) \neq 0$ for $k = 1, \cdots, l$, and let $\dot{t}_1, \dot{t}_2 \in \tilde{T}$ be any representatives of $t_1, t_2 \in T$, then $t_1, t_2 \in X_z$ is equivalent to $\dot{t}_1, \dot{t}_2$ being similar in $\operatorname{End}(H(\lambda_k))$ up to a constant for $k = 1, \cdots, l$. This is due to the following:

\begin{lemma}
\label{conjugacy equation lemma}
Let $A_1, A_2$ be two nonzero diagonal matrices in $\operatorname{End}(\mathbb{C}^m)$ such that $\operatorname{tr}(A_1) \neq 0$ and $\operatorname{tr}(A_1^e) \operatorname{tr}(A_2)^e = \operatorname{tr}(A_1)^e \operatorname{tr}(A_2^e$ for any $e = 1, \cdots, m$. Then there exists some constant $t \in \mathbb{C}^*$ such that $A_1$ is conjugate to $t A_2$.   
\end{lemma}
\begin{proof}
If $\operatorname{tr}(A_2) = 0$, then $\operatorname{tr}(A_2^e) = 0$ for all $e = 1, \cdots, m$. However, since $A_2$ is nonzero, we have $\operatorname{tr}(A_2) \neq 0$. Let $t = \frac{\operatorname{tr}(A_1)}{\operatorname{tr}(A_2)}$, then $\operatorname{tr}((t A_2)^e) = \operatorname{tr}((A_1)^e)$ for $e = 1, \cdots, m$. By the properties of symmetric polynomials, $A_1$ and $t A_2$ have the same set of eigenvalues counted with multiplicities. Therefore, there exists a permutation matrix $w$ such that $w A_1 w^{-1} = t A_2$.
\end{proof}

We also have an isomorphism $T \rightarrow (\mathbb{C}^*)^l$ sending $t$ to $(\alpha_1(t), \cdots, \alpha_l(t))$. The inverse map can be regarded as the restriction of $a$, and we will still use the notation $a$ to denote it.

\begin{proposition}
The following statements are equivalent:
\begin{enumerate}
    \item For $z \in (\mathbb{C}^*)^l$ such that $\sum_{1 \leq j \leq m_k} \left( \prod_{u=1}^{l} z_u^{\langle \lambda_u, ^k \mu_j + \lambda_k \rangle} \right) \neq 0$ for $k = 1, \cdots, l$, $X_z$ is the closure of a Steinberg fiber $\operatorname{St}^{-1}([a(z)])$.
    \item Suppose there are two semisimple elements $s_1, s_2$ of $\tilde{G}$, and $\operatorname{tr}_{H(\lambda_k)}(s_1) \neq 0$ for all $k = 1, \cdots, l$. Then $s_1$ and $s_2$ are conjugate in $\tilde{G}$ up to a multiple of an element in $Z(\tilde{G})$ if and only if they are conjugate in $\operatorname{GL}(H(\lambda_k))$ up to a constant for $k = 1, \cdots, l$.
\end{enumerate}
\end{proposition}

\begin{remark}
We still do not know whether the above proposition holds for other $z$ or not.
\end{remark}

The second statement is true for $A_l$-type groups even without the trace condition, and the first statement holds for $A_l$-type groups.

Thus, we raise the following question:

\begin{conjecture}
\label{Adjoint Steinberg conjecture}
Semisimple elements of $\tilde{G}$ are conjugate in $\tilde{G}$ up to a multiple of an element in $Z(\tilde{G})$ if and only if they are conjugate in $\operatorname{GL}(V(\lambda_k))$ up to a constant for $k = 1, \cdots, l$.
\end{conjecture}

This question is very similar to the conjectures raised by Steinberg at the 1966 International Congress of Mathematicians in Moscow \cite{steinberg1968classes}:

\begin{conjecture}[Steinberg]
Elements $x$ and $y$ of $G$ are conjugate in $G$ if $\rho(x)$ and $\rho(y)$ are similar for every irreducible rational representation $\rho$ of $G$.
\end{conjecture}

In \cite{steinberg1965regular}, Steinberg proves the following:

\begin{proposition}[Steinberg]
\label{conjugation proposition}
Semisimple elements $x$ and $y$ of $\tilde{G}$ are conjugate in $\tilde{G}$ if $\operatorname{tr}_{V(\lambda_k)}(x) = \operatorname{tr}_{V(\lambda_k)}(y)$ for every fundamental representation of $\tilde{G}$.
\end{proposition}

Thus, our conjecture resembles an adjoint form of the above proposition.

We can construct the family 
$$Y = \{(g, z) \in X \times \mathbb{C}^l \mid f_{z, e_1, e_2, k}(g) = 0, \forall 1 \leq e_1 < e_2 \leq \operatorname{dim}(H(\lambda_k)), 1 \leq k \leq l\} \subset X \times \mathbb{C}^l$$ 
and the projection to the second factor induces $\pi: Y \rightarrow \mathbb{C}^l$. 

Denote the image of $a$ by $Z = a(\mathbb{C}^l)$. As in Proposition \ref{prop6}, we have $\overline{\rho(T)} = \cup_{w \in W} wZw^{-1}$. Let $Y_w$ be the family defined over $wZw^{-1} \cong \mathbb{C}^l$, constructed in the same manner as $Y_{id} := Y$, by changing the choice of fundamental Weyl chamber. We can glue these families to obtain a larger family $Y := \cup_{w \in W} Y_w$ over $\overline{\rho(T)}/W$. 

\begin{problem}
We have a rational map $X \dashrightarrow \overline{\rho(T)}/W$ as shown in \cite{he2011semistable}. Is $Y$ the closure of the graph of this rational map?
\end{problem}

For convenience, we will still consider the family $\pi: Y \rightarrow \mathbb{C}^l$ in the following sections.
\subsection{Partial Proof of the Conjecture}\label{proof of conjecture}

\begin{lemma}
\label{trace dual lemma}
Let $(f,V)$ be a representation of $\tilde{G}$ such that $(f,V) \cong (f^*, V^*)$. If $s_1, s_2 \in \tilde{G}$ are conjugate in $\operatorname{GL}(V)$ up to a constant, then $\operatorname{tr}_{V}(f(s_1)) = \operatorname{tr}_{V}(f(s_2)) = 0$ or $f(s_1)$ is conjugate to $f(s_2)$ or $-f(s_2)$.
\end{lemma}

\begin{proof}
Suppose $f(s_1)$ is conjugate to $t f(s_2)$ for some $t \in \mathbb{C}^*$. We have $\operatorname{tr}_{V}(f(s_1)) = \operatorname{tr}_{V}(t f(s_2)) = t \operatorname{tr}_{V}(f(s_2))$. Also, $f^*(s_1)$ is conjugate to $t^{-1} f^*(s_2)$, and $\operatorname{tr}_{V^*}(f(s_1)) = t^{-1} \operatorname{tr}_{V^*}(f(s_2))$.

Since $V \cong V^*$, we have $\operatorname{tr}_{V} = \operatorname{tr}_{V^*}$, and thus $t^2 \operatorname{tr}_{V}(f(s_2)) = \operatorname{tr}_{V}(f(s_2))$. This gives us $t = \pm 1$ or $\operatorname{tr}_{V}(f(s_1)) = \operatorname{tr}_{V}(f(s_2)) = 0$.
\end{proof}

For example, if $\mathfrak{g}$ is of type $A_1, B_l, C_l, D_l$ (with $l$ even), $E_7, E_8, F_4, G_2$, then for all finite-dimensional $V$ of $\mathfrak{g}$, $V$ is isomorphic to its dual $V^*$ \cite[Problem~13.5]{humphreys2012introduction}.

\subsubsection{$G_2$}
\begin{proposition}
If $G$ is of type $G_2$, then semisimple elements of $\tilde{G}$ are conjugate in $\tilde{G}$ if and only if they are conjugate in $\operatorname{GL}(V(\lambda_k))$ up to a constant for $k = 1, 2$.
\end{proposition}

\begin{proof}
Let $V_1$ be the standard representation of $\tilde{G}$, and let $V_2$ be the other fundamental representation. Let $\rho_k := \rho_{\lambda_k}$. Suppose we have $t_1, t_2 \in \tilde{T}$ and $\rho_k(t_1)$ is conjugate to $a_k \rho_k(t_2)$ in $\operatorname{GL}(V_k)$ for some $a_k \in \mathbb{C}^*$. Then $a_1^7 = 1 = a_2^{14}$ since $G$ is simple. We have $\Lambda^2 V_1 = V_1 \oplus V_2$, so 
$$\operatorname{tr}_{V_1}(\rho_1(t_1)) + \operatorname{tr}_{V_2}(\rho_2(t_1)) = a_1^2 \left( \operatorname{tr}_{V_1}(\rho_1(t_2)) + \operatorname{tr}_{V_2}(\rho_2(t_2)) \right).$$

If $\operatorname{tr}_{V_1}(\rho_1(t_1)) \neq 0$, by Lemma \ref{trace dual lemma}, we have $a_1 = 1$ since the order of $a_1$ is odd. Then 
$$\operatorname{tr}_{V_1}(\rho_1(t_1)) + \operatorname{tr}_{V_2}(\rho_2(t_1)) = \left( \operatorname{tr}_{V_1}(\rho_1(t_2)) + \operatorname{tr}_{V_2}(\rho_2(t_2)) \right),$$ 
implying $\operatorname{tr}_{V_2}(\rho_2(t_1)) = \operatorname{tr}_{V_2}(\rho_2(t_2))$.

If $\operatorname{tr}_{V_1}(\rho_1(t_1)) = 0$ but $\operatorname{tr}_{V_2}(\rho_2(t_1)) \neq 0$, then we have $\operatorname{tr}_{\wedge^2 V_1}(\wedge^2 \rho_1(t_1)) \neq 0$ and $a_1^4 = 1$ by Lemma \ref{trace dual lemma}. Therefore, $a_1 = a_2 = 1$.

The remaining case is when $\operatorname{tr}_{V_1}(\rho_1(t_1)) = \operatorname{tr}_{V_2}(\rho_2(t_1)) = 0$. In all cases, we have $\operatorname{tr}_{V_k}(\rho_k(t_1)) = \operatorname{tr}_{V_k}(\rho_k(t_2))$, and by Proposition \ref{conjugation proposition}, we conclude that $t_1$ is conjugate to $t_2$.
\end{proof}

From the above proof, we observe that semisimple elements of $\tilde{G}$ are conjugate in $\tilde{G}$ if and only if they are conjugate in $\operatorname{GL}(V(\lambda_1))$ up to a constant.

\subsubsection{$B_l$}
\begin{proposition}
If $G = \operatorname{SO}_{2l+1}(\mathbb{C})$, then semisimple elements of $\tilde{G} = \operatorname{Spin}_{2l+1}(\mathbb{C})$ are conjugate in $\tilde{G}$ up to a multiple of an element in $Z(\tilde{G})$ if and only if they are conjugate in $\operatorname{GL}(V(\lambda_k))$ up to a constant for $k = 1, \cdots, l$.
\end{proposition}

\begin{proof}
Let $V_1 = \mathbb{C}^{2l+1}$ be the standard representation of $\operatorname{Spin}_{2l+1}(\mathbb{C})$. Let $V_k = \wedge^k V_1$ for $k = 2, \cdots, l-1$, and let $V_l$ be the fundamental spin representation. Then $V_k, k = 1, \cdots, l$ are the fundamental representations of $\operatorname{Spin}_{2l+1}(\mathbb{C})$. Let $\rho_k: \tilde{G} \rightarrow \operatorname{GL}(V_k)$ be the corresponding homomorphism. 

Suppose we have $t_1, t_2 \in \tilde{T}$ and $\rho_k(t_1)$ is conjugate to $a_k \rho_k(t_2)$ in $\operatorname{GL}(V_k)$ for some $a_k \in \mathbb{C}^*$ for $k = 1, \cdots, l$. The center of $\operatorname{Spin}_{2l+1}(\mathbb{C})$ is $\{e, \omega\}$. We know that $\rho_k(\omega) = 1$ if $k = 1, \cdots, l-1$ and $\rho_k(\omega) = -1$ if $k = l$. We may assume $a_k = a_1^k$ if $k = 1, \cdots, l-1$.

If $\operatorname{tr}_{V_k}(\rho_k(t_1)) = 0$ for $k = 1, \cdots, l-1$, then $\operatorname{tr}_{V_k}(\rho_k(t_1)) = \operatorname{tr}_{V_k}(\rho_k(t_2)) = 0$ for $k = 1, \cdots, l-1$. If $\operatorname{tr}_{V_l}(\rho_l(t_1)) = 0$, then $\operatorname{tr}_{V_l}(\rho_l(t_1)) = \operatorname{tr}_{V_l}(\rho_l(t_2)) = 0$. Therefore, $t_1$ is conjugate to $t_2$ by Proposition \ref{conjugation proposition}. If $\operatorname{tr}_{V_l}(\rho_l(t_1)) \neq 0$, then $a_l = 1$ or $a_l = -1$. If $a_l = 1$, then $t_1$ is conjugate to $t_2$. If $a_l = -1$, then $t_1$ is conjugate to $wt_2$.

If $\operatorname{tr}_{V_k}(\rho_k(t_1)) \neq 0$ for some $k \in \{1, \cdots, l-1\}$, then we must have $a_1^{k} = 1$ for its order is odd, and Lemma \ref{trace dual lemma} implies $a_1^{2k} = 1$. Therefore, $\operatorname{tr}_{V_k}(\rho_k(t_1)) = \operatorname{tr}_{V_k}(\rho_k(t_2))$. Repeating the discussion on $\operatorname{tr}_{V_l}(\rho_l(t_1))$ in the last paragraph, we conclude that $t_1$ is conjugate to $t_2$ or $wt_2$.
\end{proof}
\subsubsection{$C_l$}
\begin{proposition}
If $G = \operatorname{PSp}_{2l}(\mathbb{C})$, then semisimple elements of $\tilde{G} = \operatorname{Sp}_{2l}(\mathbb{C})$ are conjugate in $\tilde{G}$ up to a multiple of an element in $Z(\tilde{G})$ if and only if they are conjugate in $\operatorname{GL}(V(\lambda_k))$ up to a constant for $k = 1, \cdots, l$.
\end{proposition}

\begin{proof}
Let $V_1 = \mathbb{C}^{2l}$ be the standard representation of $\tilde{G}$. Let $V_k$ be the fundamental representations corresponding to the nodes in the Dynkin diagram in the order such that the double arrow points from the last. Let $\rho_k: \tilde{G} \rightarrow \operatorname{GL}(V_k)$ be the corresponding homomorphism. Suppose we have $t_1, t_2 \in \tilde{T}$ and $\rho_k(t_1)$ is conjugate to $a_k \rho_k(t_2)$ in $\operatorname{GL}(V_k)$ for some $a_k \in \mathbb{C}^*$ for $k = 1, \cdots, l$. The center of $\operatorname{Sp}_{2l}(\mathbb{C})$ is $\{e, \omega\}$. We know that $\rho_k(\omega) = (-1)^k$, and we have $\wedge^k V_1 = \wedge^{k-2} V_1 \oplus V_k$. Therefore, we need to show that 
$$\operatorname{tr}_{\wedge^k V_1}(\wedge^k \rho_1(t_1)) = (-1)^k \operatorname{tr}_{\wedge^k V_1}(\wedge^k \rho_1(t_2))$$ 
for all $k = 1, \cdots, l$ or that 
$$\operatorname{tr}_{\wedge^k V_1}(\wedge^k \rho_1(t_1)) = \operatorname{tr}_{\wedge^k V_1}(\wedge^k \rho_1(t_2))$$ 
for all $k = 1, \cdots, l$.

Let the order of $a_1$ be $r$, and assume $r > 1$ (if $r = 1$, then we are done). 

If $r$ is odd, then $\wedge^k \rho_1(t_1)$ being conjugate to $a_1^k \wedge^k \rho_1(t_2)$ implies that $a_1^k = 1$ or $\operatorname{tr}_{\wedge^k V_1}(\wedge^k \rho_1(t_1)) = \operatorname{tr}_{\wedge^k V_1}(\wedge^k \rho_1(t_2)) = 0$ by Lemma \ref{conjugation proposition}. Thus, for all $k = 1, \cdots, l$, we have 
$$\operatorname{tr}_{\wedge^k V_1}(\wedge^k \rho_1(t_1)) = \operatorname{tr}_{\wedge^k V_1}(\wedge^k \rho_1(t_2)).$$
Therefore, $t_1$ is conjugate to $t_2$.

If $r$ is even but $r = 2r'$ and $r'$ is odd, and if $\operatorname{tr}_{\wedge^k V_1}(\wedge^k \rho_1(t_1)) \neq 0$ for some $k \in \{1, \cdots, l\}$, then $a_1^{2k} = 1$ by Lemma \ref{conjugation proposition}. Therefore, $k$ is divisible by $r'$, and we have $a^k = (a^{r'})^{\frac{k}{r'}} = (-1)^{\frac{k}{r'}} = (-1)^k$. Consequently, we find that 
$$\operatorname{tr}_{\wedge^k V_1}(\wedge^k \rho_1(t_1)) = (-1)^k \operatorname{tr}_{\wedge^k V_1}(\wedge^k \rho_1(t_2))$$ 
for all $k = 1, \cdots, l$. Thus, $t_1$ is conjugate to $wt_2$.

If $r$ is divisible by $4$, let $r = 4r'$. We have $\operatorname{tr}_{\wedge^k V_1}(\wedge^k \rho_1(t_1)) = 0$ unless $k$ is divisible by $2r'$ by Lemma \ref{trace dual lemma}. 

We first claim that $a_2 = 1$ in this case. If $r' > 1$, then $\operatorname{tr}_{\wedge^2 V_1}(\wedge^2 \rho_1(t_1)) = \operatorname{tr}_{\wedge^2 V_1}(\wedge^2 \rho_1(t_2)) = 0$ since $2$ is not divisible by $2r'$. Therefore, we have 
$$0 = \operatorname{tr}_{\wedge^2 V_1}(\wedge^2 \rho_1(t_1)) = 1 + \operatorname{tr}_{V_2}(\rho_2(t_1))$$ 
and 
$$0 = \operatorname{tr}_{\wedge^2 V_1}(\wedge^2 \rho_1(t_2)) = 1 + \operatorname{tr}_{V_2}(\rho_2(t_2))$$ 
implying $a_2 = 1$. If $r' = 1$, we may assume $\operatorname{tr}_{\wedge^2 V_1}(\wedge^2 \rho_1(t_1)) \neq 0$. Then $a_1^4 = 1$ and $a_1^2 = -1$. We have 
$$\operatorname{tr}_{\wedge^2 V_1}(\wedge^2 \rho_1(t_1)) = a_1^2 \operatorname{tr}_{\wedge^2 V_1}(\wedge^2 \rho_1(t_2))$$ 
and thus 
$$1 + a_2 \operatorname{tr}_{V_2}(\rho_2(t_2)) = -1 - \operatorname{tr}_{V_2}(\rho_2(t_2)).$$ 
This leads to 
$$2 = -(1 + a_2) \operatorname{tr}_{V_2}(\rho_2(t_2)),$$ 
indicating that $\operatorname{tr}_{V_2}(\rho_2(t_2)) \neq 0$ and $a_2 \neq -1$. By Lemma \ref{trace dual lemma}, we conclude that $a_2^2 = 1$. Thus, we find $a_2 = 1$.

Suppose the eigenvalues of $\rho_1(t_1)$ are $a_1, \cdots, a_{2l}$ and the eigenvalues of $\rho_1(t_2)$ are $b_1, \cdots, b_{2l}$. Since $a_2 = 1$, we conclude that $\wedge^2 \rho_1(t_1)$ is conjugate to $\wedge^2 \rho_1(t_2)$, which implies that 
$$\sum_{i < j} a_i a_j = \sum_{i < j} b_i b_j.$$

For $x = (x_1, \cdots, x_n)$, let $p_k(x) = \sum_{i=1}^{n} x_i^k = x_1^k + \cdots + x_n^k$ be the $k$-th power sum. Denote by $e_k(x) = \sum_{1 \leq i_1 < i_2 < \cdots < i_k \leq n} x_{i_1} x_{i_2} \cdots x_{i_k}$ the $k$-th elementary symmetric polynomial. Let $\wedge^2 x = (x_ix_j)_{1 \leq i < j \leq n}$. Recall the Newton's identities:
$$
p_k(x) = (-1)^{k-1} k e_k(x) + \sum_{i=1}^{k-1} (-1)^{k-1+i} e_{k-i}(x) p_i(x)
$$
valid for all $n \geq 1$ and $n \geq k \geq 1$.

Writing $a = (a_1, \cdots, a_{2l}), b = (b_1, \cdots, b_{2l})$, we have $e_k(a) = \operatorname{tr}_{\wedge^k V_1}(\wedge^k \rho_1(t_1))$ and $e_k(b) = \operatorname{tr}_{\wedge^k V_1}(\wedge^k \rho_1(t_2))$. Since $\operatorname{tr}_{\wedge^k V_1}(\wedge^k \rho_1(t_1)) = 0$ unless $k$ is divisible by $2r'$, we find that $e_k(a) = e_k(b) = 0$ for $k$ odd. By Newton's identities, we have $p_k(a) = p_k(b) = 0$ for $k$ odd. Since $\wedge^2 \rho_1(t_1)$ is conjugate to $\wedge^2 \rho_1(t_2)$, we conclude that $e_k(\wedge^2 a) = e_k(\wedge^2 b)$ for all $k$. The question now reduces to proving the following lemma.
\end{proof}

\begin{lemma}
If we have two sets of complex numbers $a = (a_1, \cdots, a_n), b = (b_1, \cdots, b_n)$ such that $p_k(a) = p_k(b) = 0$ for $k$ odd and $e_k(\wedge^2 a) = e_k(\wedge^2 b)$ for all $k$, then $p_k(a) = p_k(b)$ for all $k$.
\end{lemma}

\begin{proof}
We prove by induction. Firstly, for $k = 2$:
$$p_2(a) = (p_1(a))^2 - 2e_2(a) = -2e_1(\wedge^2 a) = -2e_1(\wedge^2 b) = p_2(b).$$
Then for even $k$, we have:
$$p_k(a) = (p_{k/2}(a))^2 - p_{k/2}(\wedge^2 a) = (p_{k/2}(b))^2 - p_{k/2}(\wedge^2 b) = p_k(b).$$
\end{proof}
\subsection{The Fibers at the Vertices of the Polytope of $\overline{\rho(T)}$}\label{vertex}
In the last section, we constructed a surjective morphism $\pi: Y \rightarrow \mathbb{C}^l$ whose general fibers are isomorphic to the closures of one regular semisimple orbit and whose fiber over $(\mathbb{C}^*)^l$ is the Steinberg fibers (at least in the case of $A_l$-type). Note that $a(0) \in Z \subset \overline{\rho(T)}$ is one of the fixed points under the torus action. Let $p: Y \rightarrow X$ be the projection to the first factor. We can view $p(\pi^{-1}(0))$ as the degeneration of the family of closures of the regular semisimple conjugacy class parametrized by general points in $\mathbb{C}^l$. 

For the central fiber, $p(\pi^{-1}(0))$, we will be interested in the following, which will usually have dimension strictly larger than $\operatorname{dim} G/T$ (for example, when $G = \operatorname{PGL}_{l+1}(\mathbb{C})$ with $l \geq 3$).

\begin{lemma}
Let $V$ be a complex vector space. If $A \in \operatorname{End}(V)$ is a matrix such that $\operatorname{tr}(A^j)^{\operatorname{l.c.m}(j,k)/j} = \operatorname{tr}(A^k)^{\operatorname{l.c.m}(j,k)/k}$ for $1 \leq j < k \leq \operatorname{dim}(V)$, then $A$ has at most one nonzero eigenvalue counted with multiplicity.
\end{lemma}

\begin{proof}
If $\operatorname{tr}(A) = 0$, then $\operatorname{tr}(A^k) = \operatorname{tr}(A)^k = 0$ for all $1 \leq k \leq \operatorname{dim}(V)$, thus $A$ is nilpotent, and the statement holds.

If $\operatorname{tr}(A) \neq 0$, then the semisimple part of $A$ is conjugate to $E_{11}$ by Lemma \ref{conjugacy equation lemma}, where $E_{11}$ is the diagonal matrix with $1$ at the $(1,1)$-entry and $0$ at all other entries.
\end{proof}

\begin{theorem}
\label{prop9}
The central fiber $p(\pi^{-1}(0))$ is the union of $G$-stable pieces $S^w_J$ such that either $\operatorname{ht}(w(\alpha_k)) > \operatorname{ht}(\alpha_k)$ or $w(\lambda_k) \neq \lambda_k$ for each $k \in I \backslash J$. 
\end{theorem}

\begin{proof}
Since $p(\pi^{-1}(0)) = X_0$ is defined by 
$$\operatorname{tr}_{H(\lambda_k)}(\rho_{\lambda_k}(g)^{e_1})^{\operatorname{l.c.m}(e_1,e_2)/e_1} - \operatorname{tr}_{H(\lambda_k)}(\rho_{\lambda_k}(g)^{e_2})^{\operatorname{l.c.m}(e_1,e_2)/e_2} = 0$$ 
for $1 \leq e_1 < e_2 \leq \operatorname{dim}(H(\lambda_k)), 1 \leq k \leq l$, the set of eigenvalues contains only one nonzero eigenvalue (counted with multiplicity) for each fundamental representation. If for some $k$ we have $\operatorname{ht}(w(\alpha_k)) = \operatorname{ht}(\alpha_k)$ and $w(\lambda_k) = \lambda_k$ (noting that these two conditions imply $w(\alpha_k) = \alpha_k$), then $b(w,1) z_J$ will have at least the two eigenvectors $^k v_{m_k}, ^k v_{m_k-1}$ with two nonzero eigenvalues. If $\operatorname{ht}(w(\alpha_k)) > \operatorname{ht}(\alpha_k)$ or $w(\lambda_k) \neq \lambda_k$ for each $1 \leq k \leq l$, then obviously $S_J^w$ satisfies the equations.
\end{proof}

Let 
$$A_J = \{ w \in W^J \mid \operatorname{ht}(w(\alpha_k)) > \operatorname{ht}(\alpha_k) \text{ or } w(\lambda_k) \neq \lambda_k, \forall k \in I \backslash J \}.$$
Thus, we obtain that the central fiber is a union of $G$-stable pieces:
$$p(\pi^{-1}(0)) = \coprod_{J \subset I, w \in A^J} S_J^w.$$

We denote it as $F_G$ and call it the \textbf{central fiber} or \textbf{degenerate fiber}. 

\begin{remark}
Although we consider the equations only for the fundamental representations, for any other irreducible representation $H(\lambda)$, similar equations will add the condition $w(\lambda) \neq \lambda$ or $w(\alpha_k) \neq \alpha_k, k \in \operatorname{supp}(\lambda) \cap (I \backslash J)$ for $w \in W^J$. $F_G$ will still be in the central defined by equations like $\operatorname{tr}(A)^k = \operatorname{tr}(A^k)$ since $w(\lambda) = \lambda$ is equivalent to $w(\lambda_k) = \lambda_k$.
\end{remark}

By combinatorial argument, we have 

\begin{proposition}
\label{prop7}
Let $m_G = \operatorname{max}_{J \subset I, w \in A^J} \operatorname{dim} Z_J^w$, then $m_{A_l} = l(l+1) + \lfloor \frac{l}{3} \rfloor$.
\end{proposition}

\begin{proof}
If $l$ is divisible by $3$, we divide the Dynkin diagram into $l/3$ parts equally, with each part containing $3$ nodes. Then we need to consider the following two cases:

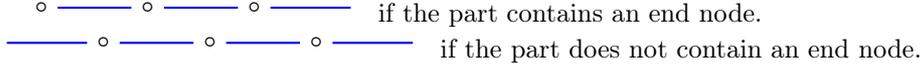
\begin{tikzpicture}[scale=.7]
  \node (1) at (-6,0) {$\circ$};
  \node (2) at (-4,0) {$\circ$};
  \node (3) at (-2,0) {$\circ$};
  \node (4) at (0,0) {};
 
  \draw [thick, blue] (1) -- (2) -- (3) -- (4);
\end{tikzpicture}
if the part contains an end node.\\
\begin{tikzpicture}[scale=.7]
  \node (0) at (-8,0) {};
  \node (1) at (-6,0) {$\circ$};
  \node (2) at (-4,0) {$\circ$};
  \node (3) at (-2,0) {$\circ$};
  \node (4) at (0,0) {};
 
  \draw [thick, blue] (0) -- (1) -- (2) -- (3) -- (4);
\end{tikzpicture}
if the part does not contain an end node.\\
We color the nodes in $J$ black; then for the first case, there are eight situations. For example, if $J$ corresponds to the following situation:
\\
\begin{tikzpicture}[scale=.7]
  \node (1) at (-6,0) {$\bullet$};
  \node (2) at (-4,0) {$\circ$};
  \node (3) at (-2,0) {$\circ$};
  \node (4) at (0,0) {};
 
  \draw [thick, blue] (1) -- (2) -- (3) -- (4);
\end{tikzpicture}
\\
Then a simple reflection corresponding to one of the nodes above must occur in the support of $w$ if $w \in A_J$, so we have at most one admissible dimension in the three connected nodes. A similar discussion applies to the other seven situations and also the second case.

If $l$ is not divisible by $3$, then we have one more case:
\\
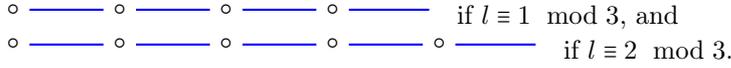
\begin{tikzpicture}[scale=.7]
  \node (1) at (-6,0) {$\circ$};
  \node (2) at (-4,0) {$\circ$};
  \node (3) at (-2,0) {$\circ$};
  \node (4) at (0,0) {$\circ$};
  \node (5) at (2,0) {};
  \draw [thick, blue] (1) -- (2) -- (3) -- (4) -- (5);
\end{tikzpicture}
if $l \equiv 1 \mod 3$, and\\
\begin{tikzpicture}[scale=.7]
  \node (1) at (-6,0) {$\circ$};
  \node (2) at (-4,0) {$\circ$};
  \node (3) at (-2,0) {$\circ$};
  \node (4) at (0,0) {$\circ$};
  \node (5) at (2,0) {$\circ$};
  \node (6) at (4,0) {};
  \draw [thick, blue] (1) -- (2) -- (3) -- (4) -- (5) -- (6);
\end{tikzpicture}
if $l \equiv 2 \mod 3$.\\
Then the Dynkin diagram is again divided into $l/3$ parts, with one more kind of part we haven't studied. But a similar argument applies, and we still have at most one admissible dimension in the four (or five) connected nodes containing an end node.
\end{proof}

If $l$ is divisible by $3$, let $J = \{2, 5, \cdots, l-1\}$. Let $w = s_2 s_5 \cdots s_{l-1}$. Then it is easy to check that $w \in A_J$ and $\operatorname{dim} S^w_J$ will reach the maximal number $m_{A_l}$. Moreover, it is unique.

\begin{proposition}
If $l$ is divisible by $3$, then $S^w_J$ for $J = \{2, 5, \cdots, l-1\}$ and $w = s_2 s_5 \cdots s_{l-1}$ is the unique $G$-stable piece in $F_{A_l}$ with the maximal dimension $m_{A_l}$.
\end{proposition}

\begin{proof}
From the proof of Proposition \ref{prop7}, we see that in each part, the admissible dimension is at most one, occurring exactly in the following situations:
\\
\begin{tikzpicture}[scale=.7]
  \node (1) at (-6,0) {$\bullet$};
  \node (2) at (-4,0) {$\circ$};
  \node (3) at (-2,0) {$\hat{\circ}$};
  \node (4) at (0,0) {};
  \draw [thick, blue] (1) -- (2) -- (3) -- (4);
\end{tikzpicture}\\
Here, a black dot indicates that $J$ contains this node, and a hat on the node indicates that the corresponding simple reflection occurs in $w$ once.

We see that if in a part there is a hat on a white circle, then there is a hat on some white circle that is connected to another part. However, none of the parts in the above form could connect to it. So a hat must lie on a black node, and such a situation is unique.
\end{proof}
\begin{proposition}
If $l \equiv 1 \operatorname{mod} 3$, then the number of $G$-stable pieces in $F_{A_l}$ with the maximal dimension $m_{A_l}$ is $l$.
\end{proposition}

\begin{proof}
The only patterns that could appear in a part of four connected nodes are the following:
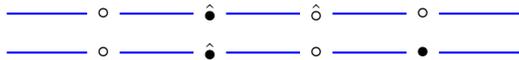
\begin{tikzpicture}[scale=.7]
    \node (0) at (-8,0) {};
    \node (1) at (-6,0) {${\circ}$};
    \node (2) at (-4,0) {$\hat{\bullet}$};
    \node (3) at (-2,0) {$\hat{\circ}$};
    \node (4) at (0,0) {$\circ$};
    \node (5) at (2,0) {};
    \draw [thick, blue] (0) --(1) -- (2) -- (3) -- (4) -- (5);
\end{tikzpicture}\\
\begin{tikzpicture}[scale=.7]
    \node (0) at (-8,0) {};
    \node (1) at (-6,0) {${\circ}$};
    \node (2) at (-4,0) {$\hat{\bullet}$};
    \node (3) at (-2,0) {${\circ}$};
    \node (4) at (0,0) {$\bullet$};
    \node (5) at (2,0) {};
    \draw [thick, blue] (0) --(1) -- (2) -- (3) -- (4) -- (5);
\end{tikzpicture}\\
Including similar patterns for the two parts containing end nodes, the first has $2\lfloor \frac{l}{3} \rfloor$ choices, and the second has $\lfloor \frac{l}{3} \rfloor + 1$ choices. Thus, there are $3\lfloor \frac{l}{3} \rfloor + 1 = l$ pieces in total.
\end{proof}

\begin{proposition}
If $l \equiv 2 \operatorname{mod} 3$, then the number of $G$-stable pieces in $F_{A_l}$ with the maximal dimension $m_{A_l}$ is $\binom{\lfloor \frac{l}{3} \rfloor + 3}{2}$.
\end{proposition}

\begin{proof}
The only extra patterns that could occur besides those when $l$ is divisible by $3$ are the following:
\\
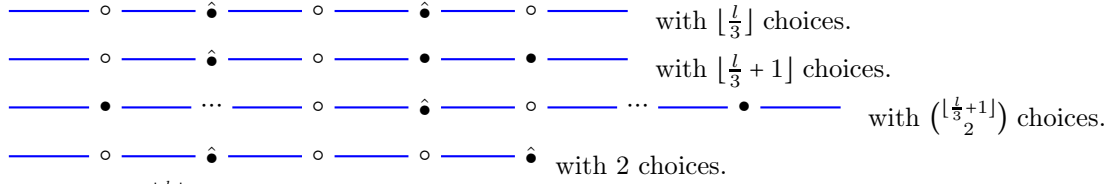
\begin{tikzpicture}[scale=.7]
    \node (0) at (-8,0) {};
    \node (1) at (-6,0) {${\circ}$};
    \node (2) at (-4,0) {$\hat{\bullet}$};
    \node (3) at (-2,0) {${\circ}$};
    \node (4) at (0,0) {$\hat{\bullet}$};
    \node (5) at (2,0) {$\circ$};
    \node (6) at (4,0) {};
    \draw [thick, blue] (0) --(1) -- (2) -- (3) -- (4) -- (5) -- (6);
\end{tikzpicture} with $\lfloor \frac{l}{3} \rfloor$ choices.\\
\begin{tikzpicture}[scale=.7]
    \node (0) at (-8,0) {};
    \node (1) at (-6,0) {${\circ}$};
    \node (2) at (-4,0) {$\hat{\bullet}$};
    \node (3) at (-2,0) {${\circ}$};
    \node (4) at (0,0) {$\bullet$};
    \node (5) at (2,0) {$\bullet$};
    \node (6) at (4,0) {};
    \draw [thick, blue] (0) --(1) -- (2) -- (3) -- (4) -- (5) -- (6);
\end{tikzpicture} with $\lfloor \frac{l}{3} + 1 \rfloor$ choices.\\
\begin{tikzpicture}[scale=.7]
    \node (0) at (-8,0) {};
    \node (1) at (-6,0) {${\bullet}$};
    \node (2) at (-4,0) {${\cdots}$};
    \node (3) at (-2,0) {${\circ}$};
    \node (4) at (0,0) {$\hat{\bullet}$};
    \node (5) at (2,0) {$\circ$};
    \node (6) at (4,0) {$\cdots$};
    \node (7) at (6,0) {$\bullet$};
    \node (8) at (8,0) {};
    \draw [thick, blue] (0) --(1) -- (2) -- (3) -- (4) -- (5) -- (6) -- (7) -- (8);
\end{tikzpicture} with $\binom{\lfloor \frac{l}{3} + 1 \rfloor}{2}$ choices.\\
\begin{tikzpicture}[scale=.7]
    \node (0) at (-8,0) {};
    \node (1) at (-6,0) {${\circ}$};
    \node (2) at (-4,0) {$\hat{\bullet}$};
    \node (3) at (-2,0) {${\circ}$};
    \node (4) at (0,0) {$\circ$};
    \node (5) at (2,0) {$\hat{\bullet}$};
    \draw [thick, blue] (0) --(1) -- (2) -- (3) -- (4) -- (5);
\end{tikzpicture} with 2 choices.\\
So there are $\binom{\lfloor \frac{l}{3} \rfloor + 3}{2}$ pieces in total.
\end{proof}

\begin{theorem}
\label{prop10}
For each $k \in I$, there exists at least one element $w_k \in A^{\{k\}}$ such that the closure of $S_{\{k\}}^{w_k}$ is an irreducible component of $F_G$.
\end{theorem}

So $F_G$ has at least $l$ irreducible components.

\begin{proof}
Choose $w_k$ to be the element $s_{j_1} s_{j_{u}} s_k$, where $j_1, \cdots, j_u$ range over the nodes that are not $k$ or the end nodes of the Dynkin diagram that are connected to node $k$ without repetition. The reflection $s_{j_s}$ is left of $s_{j_t}$ if $j_s$ and $j_t$ lie on the same interval between $k$ and some end node, and $j_s$ is farther than $j_t$ from $k$. If $F_G \supset \overline{S_J^w} \supset S_{\{k\}}^{w_k}$, then $J \subset \{k\}$ and so $J = \{k\}$. It is easy to see that $w_k$ is a subexpression of any other element in $A_{\{k\}}$, thus it is a maximal element in $A_{\{k\}}$. Hence, $S_{\{k\}}^{w_k}$ is not in the boundary of any other $G$-stable piece in $F_G$, and so the closure is an irreducible component.
\end{proof}

Let $a_\Gamma$ be the number of end nodes of the irreducible Dynkin diagram $\Gamma$. We see that $\operatorname{dim} S_{\{k\}}^{w_k} = \operatorname{dim} G/T + a_\Gamma$ if $k$ is not an end node, and $\operatorname{dim} S_{\{k\}}^{w_k} = \operatorname{dim} G/T + a_\Gamma - 1$ if $k$ is for the $w_k$ constructed in the proof of Proposition \ref{prop10}. From the proof, we also see that $F_G \cap S_{\{k\}}^\circ = \overline{S_{\{k\}}^{w_k}} \cap S_{\{k\}}^\circ$.

By definition, a variety is said to have pure dimension (or be pure dimensional) if and only if every irreducible component has the same dimension. 

\begin{corollary}\label{pure dimensional corollary}
If $G$ is simple with rank $l \geq 3$, then the degenerate fiber $F_G$ is not pure dimensional. 
\end{corollary}
\subsection{Further Discussions}\label{conjecture}

Let \( A_J^\circ := \{ w \in W^J \mid l(w) = |I \backslash J|, \text{ and } \operatorname{ht}(w(\alpha_k)) > \operatorname{ht}(\alpha_k) \text{ or } w(\lambda_k) \neq \lambda_k, \forall k \in I \backslash J \} \).

In the following, we denote \( i_G = \sum_{J \subset I} |A_J^\circ| \). If the Dynkin diagram of \( G \) is \( \Gamma \), we will also denote \( i_G \) by \( i_\Gamma \).

For example, by simple calculations, we get \( i_{A_2} = i_{B_2} = i_{G_2} = 3 \) and \( i_{A_1 \times A_1} = 1 \).

By our knowledge of symmetric groups, we could devise an algorithm to calculate \( i_{A_l} \) for a specific \( l \). The sequence \( i_{A_l} \) for \( l \leq 8 \) is \( 1, 3, 9, 26, 77, 232, 716, 2250 \).

Some problems are proposed according to the above, and we have not yet obtained satisfying answers.

\begin{problem}
The growth rate of \( i_{A_l} \) is \( 3^{l} \). Find a general formula for \( i_{A_l} \).
\end{problem}
\begin{problem}
\( |A_J^\circ| > 0 \) if \( J \) is a proper subset of \( I \).
\end{problem}

\begin{problem}
Determine and describe those regular diagonal orbits in the degenerate fiber. Does every \( G \)-stable piece \( S_J^w \) for \( w \in A_J \) contain a regular diagonal orbit?
\end{problem}

There may be other subsets in \( \mathbb{C}^l \) whose pre-image under \( \pi \) could be written as the disjoint union of \( G \)-stable pieces. 

\begin{proposition}\label{preimage of subset of base}
Let \( \Sigma_K := \{(z_1, \cdots, z_l) \mid z_i = 0 \text{ if and only if } i \in K\} \) be a \( T \)-orbit in \( \mathbb{C}^l \). Let \( B^K_J \) be the subset of \( W^J \) satisfying the following:
\begin{enumerate}
    \item For \( k \in K \backslash J \), we have \( w(\lambda_k) \neq \lambda_k \text{ or } w(\alpha_k) \neq \alpha_k \).
    \item For \( k \in J \backslash K \), we have \( w(\lambda_k) \neq \lambda_k \).
    \item For \( k \in I \backslash (J \cup K) \):
    \begin{enumerate}
        \item If \( J \cap \operatorname{comp}_{I \backslash K}(k) \neq \varnothing \), then \( w(\lambda_k) \neq \lambda_k \).
        \item If \( J \cap \operatorname{comp}_{I \backslash K}(k) = \varnothing \) and \( w(\lambda_k) = \lambda_k \), then \( w(\alpha_i) = \alpha_i \) for all \( i \in \operatorname{comp}_{I \backslash K}(k) \) and \( w(\alpha_j) \neq \alpha_j \) for all \( j \in K \backslash J \) that are connected to \( \operatorname{comp}_{I \backslash K}(k) \).
    \end{enumerate}
\end{enumerate}
Then
$$ p(\pi^{-1}(\Sigma_K)) \supset \coprod_{J \subset I, w \in B^K_J} S_J^w $$
\end{proposition}

By definition, \( \operatorname{comp}_{I \backslash K}(k) \) is the set of nodes in the subgraph generated by \( I \backslash K \) that are connected to node \( k \). In other words, it is the connected component of \( I \backslash K \) that contains \( k \). Node \( j \) is connected to \( \operatorname{comp}_{I \backslash K}(k) \) if and only if there is an edge between \( j \) and some node in \( \operatorname{comp}_{I \backslash K}(k) \) in the subgraph \( I \backslash J \).

\begin{proof}
If \( wz_J \in p(\pi^{-1}(\Sigma_K)) \). If \( k \in K \backslash J \), then the condition
$$ 
\begin{aligned}
    \operatorname{tr}_{H(\lambda_k)}(\rho_{\lambda_k}(wz_J)^{e_1})^{\operatorname{l.c.m}(e_1,e_2)/e_1} - \operatorname{tr}_{H(\lambda_k)}(\rho_{\lambda_k}(wz_J)^{e_2})^{\operatorname{l.c.m}(e_1,e_2)/e_2} = 0,\\  1 \leq e_1 < e_2 \leq \operatorname{dim}(H(\lambda_k)), \; 1 \leq k \leq l, \; \forall z \in \Sigma_K
\end{aligned}
 $$
is equivalent to \( \rho_{\lambda_k}(wz_J) \) having only one eigenspace of dimension one with a nonzero eigenvalue. This is equivalent to \( w(\lambda_k) \neq \lambda_k \text{ or } w(\alpha_k) \neq \alpha_k \).

If \( wz_J \in p(\pi^{-1}(\Sigma_K)) \) and \( k \in J \backslash K \), then \( \rho_{\lambda_k}(wz_J) \) is nilpotent or has only one eigenspace of dimension one with a nonzero eigenvalue. However, the trace conditions \( f_{z,e_1,e_2,k}(wz_J) = 0 \) require that if \( \rho_{\lambda_k}(wz_J) \) has a nonzero eigenvalue, it must have at least two eigenspaces of dimension one with nonzero eigenvalues. Therefore, \( \rho_{\lambda_k}(wz_J) \) must be nilpotent, leading to \( w(\lambda_k) \neq \lambda_k \).

For \( k \in I \backslash (J \cup K) \), if \( J \cap \operatorname{comp}_{I \backslash K}(k) \neq \varnothing \) and \( w(\lambda_k) \neq \lambda_k \), then \( \rho_{\lambda_k}(wz_J) \) is nilpotent and certainly satisfies all equations. If \( J \cap \operatorname{comp}_{I \backslash K}(k) = \varnothing \) and \( w(\lambda_k) = \lambda_k \) with \( w(\alpha_i) = \alpha_i \) for all \( i \in \operatorname{comp}_{I \backslash K}(k) \) and \( w(\alpha_j) \neq \alpha_j \) for all \( j \in K \backslash J \) that are connected to \( \operatorname{comp}_{I \backslash K}(k) \), then the eigenspaces of \( \rho_{\lambda_k}(wz_J) \) with nonzero eigenvalues are the same as those of \( \rho_{\lambda_k}(z_K) \). Since \( f_{z,e_1,e_2,k}(z_K) = 0 \) for some \( z \), it follows that \( f_{z',e_1,e_2,k}(wz_J) = 0 \) for some \( z' \).

Since \( wz_J \) is upper triangular in each fundamental representation, we see that if \( w \in B_J^K \) and some \( z \in \Sigma_K \), we have \( f_{z,e_1,e_2,k}(wz_J) = 0 \) implying \( f_{tz,e_1,e_2,k}(twz_J) = 0, \forall t \in T \). Thus, \( wz_J \in p(\pi^{-1}(\Sigma_K)) \) implies \( S_J^w \subset p(\pi^{-1}(\Sigma_K)) \) if \( w \in B_J^K \).
\end{proof}

\begin{conjecture}
$$ p(\pi^{-1}(\Sigma_K)) = \coprod_{J \subset I, w \in B^K_J} S_J^w $$
\end{conjecture}

Particularly, we know that this holds when \( K = I \).
\subsubsection{Semistable Points in the Wonderful Compactifications}\label{semistable}
The set of semistable points \( X^{ss} \) of \( X \) under the diagonal action, with the pull-back of the tautological line bundle via \( \rho_{\lambda}: X \rightarrow \mathbb{P}(\operatorname{End}(H(\lambda))) \), is studied in \cite{he2011semistable} and has a nice description as the union of the \( G \)-stable pieces:
$$ X^{ss} = \coprod_{J \subset I} S_J^e. $$
The set of semistable points \( X^{ss} \) is independent of the choice of the dominant integral weight \( \lambda \).

The quotient \( X / \!\!/ G \) is isomorphic to \( \bar{T} / W \) (\cite[Theorem~0.7]{he2011semistable}). Thus, we have a rational map \( \phi: X \dashrightarrow \bar{T}/W \) defined on \( \coprod_{J \subset I} S_J^e \).

A diagonal orbit is closed in \( X^{ss} \) if and only if it intersects \( \bar{T} \) (\cite[Lemma~0.4]{he2011semistable}). Every diagonal orbit contains a unique closed diagonal orbit, and two diagonal orbits \( O_1 \) and \( O_2 \) are in the same fiber of \( \phi \) if and only if the associated closed diagonal orbits are the same. Therefore, \( \phi \) extends the Steinberg map \( \operatorname{St}: G \rightarrow T/W \). For \( x \in T/W \), \( \phi^{-1}(x) \) is a Steinberg fiber, and for \( x \in \bar{T}/W \backslash T/W \), \( \phi^{-1}(x) \subset X \backslash G \).

The \( T \)-fixed points of \( \bar{T} \) are in an orbit under the action of \( W \). The fiber of \( \phi \) at the corresponding point in \( \bar{T}/W \) is \( S_I^e \), which equals the open Schubert cell in \( S_{i}^\circ \cong G/B \times G/B^- \). 

\subsubsection{Compactifying the Base of Another Family}\label{base compactifying}
In \cite{he2006closures}, He discusses the compactification of \( \tilde{G} \). We may denote it as \( \tilde{X} \). It is the closure of the image of the product map
$$\left(\pi, \prod_{k \in I} (\rho_{k},1)\right): \tilde{G} \rightarrow X \times \prod_{k \in I} \mathbb{P}\left(\operatorname{End}\left(H(\lambda_{k}) \oplus \mathbb{C}\right)\right).$$

There is a lemma describing the difference between \( \tilde{X} \) and \( X \):

\begin{lemma}[\cite{he2006closures}]
The projective morphism \( \pi: \tilde{X} \rightarrow X \) defines a bijection between \( \tilde{X} - \tilde{G} \) and \( X - G \). In particular, \( \pi \) is a finite morphism.
\end{lemma} 

Let \( \operatorname{tr}_{k} \) denote the trace function on \( \operatorname{End}(H(\lambda_k)) \). To each \( a_{k} \in \mathbb{C} \), we may associate a global section \( \left(\operatorname{tr}_{k}, a_{k}\right) \) of the line bundle \( \mathcal{O}_{k}(1) := \mathcal{O}_{\mathbb{P}\left(\operatorname{End}\left(H(\lambda_k)\right) \oplus \mathbb{C}\right)}(1) \) on \( \mathbb{P}\left(\operatorname{End}\left(H(\lambda_k)\right) \oplus \mathbb{C}\right) \). The pull-back of \( \left(\operatorname{tr}_{k}, a_{k}\right) \) to \( \tilde{X} \) by the morphism \( (\rho_k, 1) \) is then a global section \( f_{k, a_{k}} \) of a line bundle on \( \tilde{X} \).

Let \( \tilde{Y} = \{(g,z) \mid f_{k,z_k}(g) = 0, k \in I \} \subset \tilde{X} \times \mathbb{C}^l \). Let \( \tilde{\pi}: \tilde{Y} \rightarrow \mathbb{C}^l \) be the second projection. Then each fiber of \( \tilde{\pi} \) is the closure of a Steinberg fiber of \( \tilde{G} \) (\cite{he2006closures}). To see the degeneration at infinity, we could compactify the base \( \mathbb{C}^l \) in many ways. For example, consider \( \mathbb{C}^l \subset (\mathbb{P}^1)^l = \prod_{k \in I} \{[z_k:u_k] \in \mathbb{P}^1\} \). Let \( f_{k,u_k,z_k} \) be the pull-back of \( (u_k \operatorname{tr}_k, z_k) \). The zero set of \( f_{k,u_k,z_k} \) does not depend on the choice of the representative of \( [u_k:z_k] \in \mathbb{P}^1 \). Define \( \tilde{Y}' = \{(g,[z:u]) \mid f_{k,u_k,z_k}(g) = 0, k \in I \} \subset \tilde{X} \times (\mathbb{P}^1)^l \) and let \( \tilde{\pi}': \tilde{Y}' \rightarrow (\mathbb{P}^1)^l \) be the second projection.

Note that \( (\rho_k,1)(g) = [(\rho_k(g),0)] \) if \( g \in \tilde{X} \backslash \tilde{G} \). We have the following proposition:

\begin{proposition}
Let \( p_J \in (\mathbb{P}^1)^l \) be the point where \( z_k = 0 \) if \( k \in J \) and \( u_k = 0 \) if \( k \notin J \). Then for \( J \subsetneq I \), 
$$ \tilde{\pi}'(p_J) = \{ g \in \tilde{X} \backslash \tilde{G} \mid \operatorname{tr}_k(\rho_k(g)) = 0, k \in J \}. $$
\end{proposition}

To make an analogous construction for the adjoint group \( G \), consider the representation \( H(\Lambda) \) of \( \tilde{G} \) that is also a representation of \( G \). Let \( \Lambda_1, \cdots, \Lambda_m \) be a minimal set of generators of such \( \Lambda \). Then \( \operatorname{tr}_{H(\Lambda_1)}, \cdots, \operatorname{tr}_{H(\Lambda_m)} \) forms a generating set of \( \mathbb{C}[G]^G \), and the image of 
$$(\operatorname{tr}_{H(\Lambda_1)}, \cdots, \operatorname{tr}_{H(\Lambda_m)}): G \rightarrow \mathbb{C}^m$$ 
is \( T/W \). Thus, we could construct a similar family with the sections \( f_{i,z_i}, i = 1, \cdots, m \) over \( T/W \subset \mathbb{C}^m \) whose fiber is the closure of a Steinberg fiber of \( G \) in \( X \). The compactifications of \( T/W \) in the way we address \( \tilde{G} \) would also possibly capture the degeneration at infinity.




\bibliographystyle{plain}

\begin{thebibliography}{10}

\bibitem{bate2019orbit}
Michael Bate, Haralampos Geranios, and Benjamin Martin.
\newblock Orbit closures and invariants.
\newblock {\em Mathematische Zeitschrift}, 293(3):1121--1159, 2019.

\bibitem{batyrev2011functor}
Victor Batyrev and Mark Blume.
\newblock The functor of toric varieties associated with {W}eyl chambers and {L}osev-{M}anin moduli spaces.
\newblock {\em Tohoku Mathematical Journal, Second Series}, 63(4):581--604, 2011.

\bibitem{cox2011toric}
David~A Cox, John~B Little, and Henry~K Schenck.
\newblock {\em Toric varieties}, volume 124.
\newblock American Mathematical Soc., 2011.

\bibitem{esposito2012closures}
Francesco Esposito.
\newblock Closures of orbits under the diagonal action in the wonderful compactification of {PGL}(3).
\newblock {\em Communications in Algebra}, 40(8):3127--3140, 2012.

\bibitem{evens2008wonderful}
Sam Evens and Benjamin~F Jones.
\newblock On the wonderful compactification.
\newblock {\em arXiv preprint arXiv:0801.0456}, 2008.

\bibitem{fulton2013representation}
William Fulton and Joe Harris.
\newblock {\em Representation theory: a first course}, volume 129.
\newblock Springer Science \& Business Media, 2013.

\bibitem{he2006unipotent}
Xuhua He.
\newblock Unipotent variety in the group compactification.
\newblock {\em Advances in Mathematics}, 203(1):109--131, 2006.

\bibitem{he2007G}
Xuhua He.
\newblock The {G}-stable pieces of the wonderful compactification.
\newblock {\em Transactions of the American Mathematical Society}, 359(7):3005--3024, 2007.

\bibitem{he2011semistable}
Xuhua He and Jason Starr.
\newblock Semistable locus of a group compactification.
\newblock {\em Representation Theory of the American Mathematical Society}, 15(17):574--583, 2011.

\bibitem{he2006closures}
Xuhua He and Jesper~Funch Thomsen.
\newblock Closures of {S}teinberg fibers in twisted wonderful compactifications.
\newblock {\em Transformation groups}, 11(3):427--438, 2006.

\bibitem{humphreys1995conjugacy}
James~E Humphreys.
\newblock {\em Conjugacy classes in semisimple algebraic groups}.
\newblock Number~43. American Mathematical Soc., 1995.

\bibitem{humphreys2012introduction}
James~E Humphreys.
\newblock {\em Introduction to Lie algebras and representation theory}, volume~9.
\newblock Springer Science \& Business Media, 2012.

\bibitem{kostant1963lie}
Bertram Kostant.
\newblock Lie group representations on polynomial rings.
\newblock {\em American Journal of Mathematics}, 85(3):327--404, 1963.

\bibitem{lusztig2004parabolic}
George Lusztig.
\newblock {Parabolic character sheaves. II}.
\newblock {\em Moscow Mathematical Journal}, 4(4):869--896, 2004.

\bibitem{speyer2009matroid}
David~E Speyer.
\newblock A matroid invariant via the {K}-theory of the {G}rassmannian.
\newblock {\em Advances in Mathematics}, 221(3):882--913, 2009.

\bibitem{springer2006some}
Tonny~A Springer.
\newblock Some results on compactifications of semisimple groups.
\newblock In {\em International congress of mathematicians}, volume~2, pages 1337--1348. Citeseer, 2006.

\bibitem{steinberg1965regular}
Robert Steinberg.
\newblock Regular elements of semi-simple algebraic groups.
\newblock {\em Publications Math{\'e}matiques de l'IH{\'E}S}, 25:49--80, 1965.

\bibitem{steinberg1968classes}
Robert Steinberg.
\newblock Classes of elements of semisimple algebraic groups.
\newblock {\em Proc. Internat. Congr. Math. (Moscow, 1966)}, pages 277--284, 1968.

\end{thebibliography}

 ~\newline

\end{document}